# A Musielak–Orlicz approach for modeling uncertainties in long-memory processes

Hidekazu Yoshioka[1,*]

[1] Graduate School of Advanced Science and Technology, Japan Advanced Institute of Science and Technology, 1-1 Asahidai, Nomi, Ishikawa, Japan

* Corresponding author: yoshih@jaist.ac.jp, ORCID: 0000-0002-5293-3246

***Abstract***

This paper proposes a novel mathematical framework for modeling uncertainties in supOU processes, a common model for long-memory phenomena. We address uncertainties as distortions in reversion and Lévy measures, evaluating them simultaneously via state-dependent divergence functions on Musielak–Orlicz spaces. The core of our approach involves solving optimization problems to determine the upper- and lower-bounds of cumulants under a prescribed uncertainty set. Notably, we demonstrate that while classical measures like Kullback–Leibler divergence fail in this context, Musielak–Orlicz spaces effectively resolve these issues. Along with providing sufficient conditions for the well-posedness of these optimizations, we demonstrate the framework's practical utility through a water environmental application, modeling streamflow discharge. This work offers both a theoretical advancement and a robust tool for long-memory process analysis.

***Statements & declarations***

**Acknowledgments:** The author is very grateful for the helpful comments from the two anonymous reviewers.

**Funding:** This study was supported by the Japan Society for the Promotion of Science (KAKENHI No. 25K00240, 25K07931) and the Japan Science and Technology Agency (PRESTO No. JPMJPR24KE).

**Conflict of Interests:** The authors have no relevant financial or nonfinancial interests to disclose.

**Data Availability:** The data will be made available upon reasonable request to the corresponding author.

**Declaration of Generative AI in Scientific Writing:** The author did not use generative AI to obtain the results of this study.

## 1. Introduction

We first review existing studies and fill this research gap (**Section 1.1**). Afterward, we explain the contributions of this paper to fill this gap (**Section 1.2**). Finally, we present the outline of this paper (**Section 1.3**).

### 1.1 Study background and motivation

Long-memory processes are stochastic processes with autocorrelation functions that decay at a polynomial speed [1] in contrast to short-memory processes, including typical Markovian processes. There are a wide variety of research fields in which long-memory processes play a central role, which include but are not limited to anomalous diffusion in heterogeneous media [2], early detection of the epidemic in aquaculture [3], option pricing in finance [4], quantum dots in nanocrystals [5], electricity price forecasting [6], and gradient descent for optimization and machine learning [7].

A nominal model for long-memory phenomena is the superposition of Ornstein–Uhlenbeck processes, the so-called supOU process [8,9], which is formally a summation or integration of mutually independent Ornstein–Uhlenbeck (OU) processes with distinct reversion speeds with respect to Lévy bases[1] [10]. OU processes are representative stochastic differential equations (SDEs) whose properties have been studied extensively [e.g., 11-13], and this fact has facilitated the analysis of the supOU process. Despite its simplicity, supOU processes have rich mathematical structures because of their diverse long-memory and fluctuation properties encoded in Lévy bases. To date, supOU processes have been applied to volatilities in economics [14] and environmental [15] and ecological problems [16]. Their memory and spectral characteristics have been investigated by Kovtun et al. [17], and growth and extremes have been investigated by Grahovac and Kevei [18] and Moser and Stelzer [19]. Leonenko and Pepelyshev [20] discussed theoretical and computational methods to efficiently simulate jump-driven supOU processes with prescribed probability density functions (PDFs). Long-memory processes related to supOU processes, such as the Volterra and Trawl processes [21] as well as the superpositions of nonlinear SDEs [22,23], have also been theoretically studied.

This paper focuses on jump-driven supOU processes because they frequently appear in applied problems where the nonnegativity of solutions is necessary, such as the case for volatility and streamflow discharge [14-16] and environmental problems [23], and because the framework of uncertainty modeling in this paper straightforwardly applies to cases with diffusion. The behavior and mathematical properties of a jump-driven supOU process are determined by two measures: the reversion measure to distribute reversion speeds and the Lévy measure that governs the size and frequency of jumps. The singularity of a reversion measure determines the memory strength, namely, the decay of the autocorrelation function, of a supOU process, whereas the Lévy measure quantifies the moments and path regularities of this process. Therefore, perturbing these measures results in qualitatively and/or quantitatively different supOU processes in theory and different statistical predictions in application. Because the estimation of reversion and Lévy measures amounts to the identification of a supOU process, quantifying the relationship between uncertainties (i.e., modeling errors) introduced in these measures and the resulting model is important in both theory and application.

[1] Lévy bases are space-time generalization of Lévy processes and are seen as high-dimensional Poison random measures for pure-jump cases [9].

Uncertainties in stochastic models may come from data and assumptions. First, the available data in an applied problem contain measurement errors due to technical restrictions [24,25]. Second, the quantity of data, such as the length of time series, is usually restricted [26,27]. Third, the assumed model may deviate from the ground truth because improper functional forms are assumed [28,29]. Finally, even if the data are complete in quality and quantity, the identification procedure may result in source error. In any case, quantifying and, if possible, predicting the possible range of uncertainties is a key step toward operating a misspecified stochastic process model.

Thus far, modeling uncertainties in reversion and Lévy measures have been carried out separately, assuming that one of them is correct [30,31]; however, this assumption is too strong because both measures would be estimated on the basis of a common time series dataset for the target phenomenon. This means that there is a gap between the theory and application for modeling uncertainties in supOU processes, motivating us to develop a unified mathematical framework that can simultaneously address uncertainties contained in reversion and Lévy measures.

### 1.2 Aim and contribution

The aims of this paper are to formulate and analyze a unified framework for addressing uncertainties in jump-driven supOU processes and present its applications, with which we will be able to abandon the assumption that uncertainties are concentrated solely on either a reversion measure or a Lévy measure. An important advantage of the proposed mathematical framework is its ability to evaluate the influences of uncertainties on the statistics of supOU processes, which are non-Markovian processes, without resorting to Monte Carlo simulations.

We focus on the state-dependent uncertainty assessment approach recently proposed in Strati [32], whose primary direction was insurance, while their mathematical tools, especially Musielak–Orlicz spaces (Chapter 3 in Chlebicka et al. [33]), carry over to any measures having finite masses; reversion and Lévy measures of supOU processes are no exception. We exploit the useful features of supOU processes such that its cumulants and autocorrelations are given by certain (non)polynomial moments of the two measures, which provide the relationship to clearly show how uncertainties in these measures propagate to the statistics of supOU processes.

A meaningful method for estimating the impacts of uncertainties on a statistic is to evaluate the size of uncertainties by using a divergence, also called statistical divergence, which quantifies the difference between two different stochastic models [34]. This methodology was reviewed in Ben-Tal and Teboulle [35] and was later systematized in Fröhlich and Williamson [36] on the basis of Orlicz spaces, which are function spaces of random variables that have a certain strength of singularity, such as nonpolynomial integrability (Chapter 13 in Rubshtein et al. [37]). Finding a proper Orlicz space to which the target random variable belongs allows for the construction of divergences that are well defined and computable. We focus on the state-dependence in divergences, but qualitatively similar approaches would apply to related risk indices such as robust return risk measures [38,39] and Orlicz premia [40,41], where the well-posedness of optimization problems is determined by Orlicz spaces.

The limitation of the approach based on classical Orlicz spaces is that it cannot address state-dependent divergences. Here, the state-dependence indicates the situation where a parameter or a coefficient of divergence depends on the randomness coming from other than the target random variable. For the case of supOU processes, this corresponds to the situation in which uncertainties affect both reversion and Lévy measures and in which accounting for the state-dependence is essential because of their cumulant forms. Strati [32] proposed using Musielak–Orlicz spaces instead of classical Orlicz spaces, with which the state-dependence can be naturally formulated. We employ this recently developed theory to formulate both upper and lower bounds of cumulants of supOU processes, with which we can estimate the cumulants subject to uncertainties by finding proper solutions to optimization problems formulated with Musielak–Orlicz spaces. We show that the classical divergences, such as the Kullback–Leibler (KL) divergence, fail to define these optimization problems well and that allowing for its state-dependence resolves this issue. We also discuss what state-dependence is allowed for the alpha divergence: a popular generalization of the KL one [e.g., 42,43]. In summary, the novelty of the theory in this paper is modeling uncertainties in supOU processes by accounting for state dependence that was ignored in previous studies [30,31], which is achieved by the use of suitable Musielak–Orlicz spaces.

In our framework, the supOU process under uncertainty, i.e., the distorted supOU process, differs from the classical processes because the reversion and Lévy measures are interdependent in the former but not in the latter. This implies that we need to check the integrability condition to verify the well-posedness of ambit fields (Proposition 39 in Barndorff-Nielsen et al. [21]), and moreover, its connection to the solvability of our optimization problems should be clarified. We show how the conditions for well-posedness are linked between the distorted supOU process and our optimization problems by focusing on state-dependence.

We also contribute to the computational aspects of modeling uncertainties in supOU processes. Our optimization problems to evaluate the cumulants subject to uncertainties can be computationally solved by a gradient descent method with momentum along with a quantization method to discretize expectations. This method provides not only the solution but also the optimal Lagrangian multipliers (Chapter 5.6.2 in Boyd and Vandenberghe [44]), which are sensitive to constraints about the uncertainty set and hence are useful for a deeper understanding of the influences of uncertainties in the proposed framework.

Finally, we present demonstrative applications of the proposed mathematical framework to supOU processes fitted against real streamflow discharge data. Indeed, streamflow discharge, as a physical quantity characterizing water quantity in river environments, is characterized by polynomial and hence slow recessions after flood events [45-47], necessitating the use of a long-memory process for mathematical modeling. Moreover, each flood event can be considered a jump event [48-50]. We analyze how the reversion and Lévy measures are distorted and subject to uncertainties and how the resulting supOU processes behave in terms of both statistics and threshold dynamics. Consequently, this paper contributes to the formulation through the application of uncertainty modeling for long-memory processes by focusing on supOU processes.

### 1.3 Organization of this paper

**Section 2** introduces the supOU process and its fundamental mathematical results about its existence. **Section 3** presents our mathematical framework for modeling uncertainties in supOU processes and theoretically analyzes a pair of optimization problems based on Musielak–Orlicz spaces. **Section 4** applies our mathematical framework to real data of streamflow discharge. **Section 5** concludes this paper and presents perspectives on our study. **Appendix** presents proofs of propositions and supplementary results.

## 2. supOU process

### 2.1 Formulation

In this subsection, we formulate a supOU process whose Lévy base admits a nonclassical structure where the influences of reversion and jump are not decomposable. We work on a filtered complete probability space as usual in stochastic calculus for continuous-time processes [e.g., 14].

Time $t$ is a real independent variable. As a supOU process, we consider the following continuous-time stochastic process $X = (X_t)_{t \in \mathbb{R}}$ driven by jumps:

$$X_t = \int_{r=0}^{r=+\infty} \int_{z=0}^{z=+\infty} \int_{s=-\infty}^{s=+\infty} \underbrace{z}_{\text{Jump size}} \underbrace{e^{-r(t-s)} \mathbb{I}(s < t)}_{\text{Reversion}} \underbrace{N_\phi(\mathrm{d}r\mathrm{d}z\mathrm{d}s)}_{\text{Randomness}}, \quad t \in \mathbb{R}, \tag{1}$$

where $\mathbb{I}(s<t)$ is the indicator function such that $\mathbb{I}(s<t)=1$ if $s<t$ and $\mathbb{I}(s<t)=0$ otherwise, $r>0$ represents distributed reversion speeds, $z>0$ represents distributed jump sizes, $N$ represents a Poisson random measure on $(0,+\infty)\times(0,+\infty)\times\mathbb{R}$ with the compensated measure $\phi(r,z)\pi(\mathrm{d}r)\nu(\mathrm{d}z)\mathrm{d}s$, with $\pi$ being a probability measure of a positive random variable and $\nu(\mathrm{d}z)$ being a Lévy measure satisfying the integrability condition $\int_0^{+\infty} z\nu(\mathrm{d}z) < +\infty$ (necessary for the existence of the average of Eq.(1) when $\phi \equiv 1$), and $\phi$ is a nonnegative function on $(0,+\infty)\times(0,+\infty)$. Note that $\nu(\mathrm{d}z)$ needs not be a probability measure but $\pi(\mathrm{d}r)$ does because the latter is a probability distribution of reversion speeds, while the former is essentially a product of jump intensity and a probability distribution of jump size.

The supOU process (1) is non-Markovian because its right-hand side can be understood as an aggregation of mutually independent background OU processes [8]. More specifically, the expectation of the current value $X_t$ conditioned on the past value $X_s$ is different from that on all the background OU processes at time $s$; the information by the former is smaller than that by the latter, which gives the non-Markovian nature of the supOU process.

***Remark 1*** We recover the classical supOU process with $\phi \equiv 1$ [8].

***Remark 2*** In **Section 3**, we encounter $\phi = z^{m_1} r^{-m_2}$ and $\phi = \exp\left(m_3 z^{m_1} r^{-m_2}\right)$ with $m_1, m_2 > 0$ and $m_3 \in \mathbb{R}$ up to multiplication constants.

In classical supOU processes ( $\phi \equiv 1$ ), the compensated measure of the Poisson random measure has a separable structure between reversion and Lévy measures, while Eq.(1) does not due to the appearance of the bivariate function $\phi$. Physically, this $\phi$ implies interdependence between reversion and jump, e.g., a model where the distributions of jump size and frequency depend on the reversion speed or the distribution of reversion speeds depends on the jump size. In our context, nonconstant $\phi$ cases arise because modeling errors that simultaneously affect the reversion and jump of supOU processes are assumed. More specifically, Eq.(1) arises as a distorted model based on the assumption that a classical supOU process is identified as a model estimated from data, whereas in reality, the influences of reversion and jump may not be decomposable.

### 2.2 Mathematical structure

We need to ensure under what condition the right-hand side of Eq.(1) exists because we are not aware of literature that explicitly addresses supOU processes of this form. We exploit the integrability condition for ambit fields and adapt it to Eq.(1) (Proposition 39 in Barndorff-Nielsen et al. [21]).

***Proposition 1***

*The process* $X$ *in Eq.(1) exists and is stationary if*

$$\int_{r=0}^{r=+\infty} \int_{z=0}^{z=+\infty} \int_{s=0}^{s=+\infty} \min\left\{1, z^2 e^{-2rs}\right\} \phi(r,z) \pi(\mathrm{d}r) \nu(\mathrm{d}z) \mathrm{d}s < +\infty . \tag{2}$$

*Moreover, the characteristic function* $\mathbb{E}\left[e^{\iota u X_t}\right]$ *(* $t \in \mathbb{R}$ *) is independent of* $t$ *and is given by (* $\iota$ *represents the imaginary unit:* $\iota^2 = -1$*)*

$$\ln \mathbb{E}\left[e^{\iota u X_t}\right] = \int_{r=0}^{r=+\infty} \int_{z=0}^{z=+\infty} \int_{s=0}^{s=+\infty} \left(\exp\left(\iota u e^{-rs} z\right) - 1\right) \phi(r,z) \pi(\mathrm{d}r) \nu(\mathrm{d}z) \mathrm{d}s , \quad u \in \mathbb{R} . \tag{3}$$

For the classical case ( $\phi \equiv 1$ ), the integrability condition (2) is satisfied if the inverse moment of $\pi$ exists:

***Assumption 1***

$$\int_0^{+\infty} r^{-1} \pi(\mathrm{d}r) < +\infty . \tag{4}$$

Under **Assumption 1**, the classical supOU process is stationary and infinitely divisible (Theorem 3.1 in Barndorff-Nielsen and Stelzer [9]). Another key assumption that we employ is about densities:

***Assumption 2*** *Measures* $\pi$ *and* $\nu$ *admit positive densities.*

This is a technical assumption to ensure the existence of solutions to our optimization problems but is satisfied by many examples in the literature [e.g., 20,30]. We assume that **Assumptions 1 and 2** hold true in the rest of this paper.

For later use, we set the following integrability conditions:

**(I$_k$)** $$I_k = \frac{1}{k}\int_{r=0}^{r=+\infty}\int_{z=0}^{z=+\infty}\frac{\phi(r,z)z^k}{r}\pi(\mathrm{d}r)\nu(\mathrm{d}z) < +\infty,\quad k\in\mathbb{N}. \tag{5}$$

The condition (2) is satisfied if **(I$_2$)** is, because

$$\int_{r=0}^{r=+\infty}\int_{z=0}^{z=+\infty}\int_{s=0}^{s=+\infty}\min\left\{1,z^2e^{-2rs}\right\}\phi(r,z)\pi(\mathrm{d}r)\nu(\mathrm{d}z)\mathrm{d}s \le \int_{r=0}^{r=+\infty}\int_{z=0}^{z=+\infty}\frac{\phi(r,z)z^2}{2r}\pi(\mathrm{d}r)\nu(\mathrm{d}z). \tag{6}$$

The cumulants and autocorrelation function of Eq.(1) are derived explicitly. For example, because the left-hand side of Eq.(3) represents a cumulant generating function, we have

$$\mathrm{Cum}_k = I_k \quad \text{under } \textbf{(I}_k\textbf{)},\quad k\in\mathbb{N} \tag{7}$$

and

$$\mathrm{ACF}(h) = \frac{\mathbb{E}\left[\left(X_t-\mathbb{E}[X_t]\right)\left(X_{t+h}-\mathbb{E}[X_t]\right)\right]}{\mathrm{Cum}_2} = I_2^{-1}\int_{r=0}^{r=+\infty}\int_{z=0}^{z=+\infty}\frac{\phi(r,z)z^2}{2r}e^{-rh}\pi(\mathrm{d}r)\nu(\mathrm{d}z) \quad \text{under } \textbf{(I}_2\textbf{)}, \tag{8}$$

where $\mathrm{Cum}_k$ represents the $k$ th order cumulant and $\mathrm{ACF}(h)$ is the autocorrelation function with lag $h \ge 0$. In particular, the cases $k = 1,2,3$ represent average, variance, and nonnormalized skewness, respectively.

By a Tauberian argument [51], we can evaluate the tail of the autocorrelation function according to singularities of $\phi$ and $\pi$ when $h$ is large, as shown in the next proposition. The asymptotic condition (9) means that Eq.(1), as for the classical case, has a long-memory property if the integrand in Eq.(8) has a certain range of singularities near $r = 0$. The difference is that the singularity depends on both $\phi$ and $\pi$ in the present case.

***Proposition 2***

*Assume the condition **(I$_2$)** and that $\pi$ admits a density $p$. If the following asymptotic estimate holds true for some constant $A > 1$:*

$$\frac{\phi(r,z)p(r)}{r} \sim O\left(r^{A-2}\right) \quad \textit{with respect to } 0 < r << 1 \textit{ and all } z > 0, \tag{9}$$

*then*

$$\mathrm{ACF}(h) = O\left(h^{-(A-1)}\right) \quad \textit{for } h >> 1. \tag{10}$$

In particular, this proposition implies that the tail behavior of the autocorrelation function is the same for all $\phi$ that are bounded near $r = 0$. Note that $A$ serves as the Hurst index if $A \in (1,2]$.

***Remark 3*** According to Eqs.(5), (7), (8), and **Propositions 1-2**, the condition **(I₂)** is sufficient for the stationarity of the supOU process (1) with bounded variance.

## 3. The Musielak–Orlicz framework

### 3.1 Problem setting

We consider a situation where we have identified a classical supOU process (Eq.(1) with $\phi \equiv 1$, called the benchmark model in the sequel) from a limited dataset, but the true model corresponds to $\phi \neq 1$ (distorted model). In this case, we have estimations of the reversion measure $\pi$ and Lévy measure $\nu$ but are not necessarily correct, and the difference between the benchmark and distorted models is evaluated by $\phi$.

#### 3.1.1 Generalized alpha divergence

In this paper, the difference between the benchmark and distorted models is evaluated through the generalized alpha divergence $\mathbb{D}(\phi)$: for positive and continuous functions $w, \alpha$ on $(0,+\infty)\times(0,+\infty)$,

$$\mathbb{D}(\phi) = \int_{r=0}^{r=+\infty}\int_{z=0}^{z=+\infty} \underbrace{\Phi\left(\underbrace{w(r,z)}_{\text{Weighting}},\underbrace{\alpha(r,z)}_{\text{Shape}},\underbrace{\phi(r,z)}_{\text{Uncertainty}}\right)}_{\text{Divergence function}} \pi(\mathrm{d}r)\nu(\mathrm{d}z) \tag{11}$$

with $\Phi$ given by

$$\Phi(w_0,\alpha_0,\phi_0) = \begin{cases} w_0 \dfrac{(\phi_0)^{\alpha_0} - \alpha_0(\phi_0 - 1) - 1}{\alpha_0(\alpha_0 - 1)} & (\alpha_0 \neq 1) \\ w_0(\phi_0 \ln\phi_0 - \phi_0 + 1) & (\alpha_0 = 1) \end{cases}, \quad w_0, \alpha_0 > 0,\ \phi_0 \geq 0. \tag{12}$$

We set $\Phi(w_0,\alpha_0,\phi_0) = +\infty$ for $\phi_0 < 0$, which extends the domain of $\Phi$ to $\phi_0 \in \mathbb{R}$ without losing convexity. This extension plays a role in the derivation of the dual representations of our optimization problems. The divergence $\mathbb{D}(\phi)$ is used to evaluate the difference between the benchmark and distorted models; $\Phi = 0$ if there is no uncertainty ($\phi \equiv 1$) and is positive otherwise. The two coefficients $\alpha$ and $w$ play distinct roles; the former modulates the shape of divergence (sensitivity to uncertainties), whereas the latter modulates singularity (more connected to the well-posedness of the proposed framework), as theoretically and computationally discussed in **Sections 3.2 and 4.3**.

Divergence is a natural measure of uncertainty, whereas classical uncertainties that have no state dependence possibly fail to model the uncertainties of supOU processes. Therefore, with Eq.(11), we generalize the classical alpha divergence (e.g., Eq. (1) in Póczos and Schneider [43]) in three ways. First, we extend the domain of divergence from a PDF to a positive measurable function because of our modeling assumption that the function $\phi$ in Eq.(1) needs not be a PDF. Second, we assume that the parameter $\alpha$ depends on $r, z$, while the classical alpha divergence assumes that $\alpha$ is a positive constant, which offers a more flexible framework for modeling uncertainties. As shown in Eq.(12), the generalized alpha divergence is reduced to the classical KL one when $\alpha \equiv 1$ and $w \equiv 1$. Finally, the function $w$ serving as a weighting factor further increases the flexibility in modeling uncertainties, but its true value is revealed with the modified KL divergence function ($\alpha \equiv 1$ but $w \neq 1$), with which our optimization problems become well defined; the classical KL case never succeeds in the supOU process.

We summarize the properties of $\Phi$ in **Lemma 1**, which follow from direct calculations:

***Lemma 1***

*For each* $w_0, \alpha_0 > 0$, $\Phi(w_0, \alpha_0, \phi_0)$ *as a function of* $\phi_0 \geq 0$ *is continuous, convex, nonnegative, and globally minimized at* $\phi_0 = 1$ *with the minimum value* $\Phi(w_0, \alpha_0, 1) = 0$. *Moreover,* $\lim_{\phi_0 \to +0} \Phi(w_0, \alpha_0, \phi_0) < +\infty$, $\lim_{\phi_0 \to +\infty} \Phi(w_0, \alpha_0, \phi_0) = +\infty$ *and* $\lim_{\phi_0 \to +\infty} \frac{\Phi(w_0, \alpha_0, \phi_0)}{\phi_0} = +\infty$ *when* $\alpha_0 \geq 1$.

### 3.1.2 Optimization problems

We consider the following pair of optimization problems where we emphasize the dependence of the cumulant $\mathrm{Cum}_k = I_k$ on $\phi$: for each $k, m = 0,1,2,...$ and $\varepsilon > 0$,

**Upper-bound case** Find $$\overline{I}_{k,m} = \sup_{\phi \in \mathfrak{A}_{m,\varepsilon}} \mathrm{Cum}_k(\phi) = \sup_{\phi \in \mathfrak{A}_{m,\varepsilon}} I_k(\phi) \quad (13)$$

**Lower-bound case** Find $$\underline{I}_{k,m} = \inf_{\phi \in \mathfrak{A}_{m,\varepsilon}} \mathrm{Cum}_k(\phi) = \inf_{\phi \in \mathfrak{A}_{m,\varepsilon}} I_k(\phi) \quad (14)$$

Here, $\mathfrak{A}_{m,\varepsilon}$ is the following admissible set of measurable functions:

$$\text{Find } \mathfrak{A}_{m,\varepsilon} = \left\{ \phi : (0,+\infty)^2 \to [0,+\infty) \middle| I_m(\phi) = c_m,\ \mathbb{D}(\phi) \leq \varepsilon \right\}, \quad (15)$$

where $c_m > 0$ is the $m$th-order cumulant in the benchmark case:

$$c_m = \frac{1}{m} \int_{r=0}^{r=+\infty} \int_{z=0}^{z=+\infty} \frac{z^m}{r} \pi(\mathrm{d}r) \nu(\mathrm{d}z). \quad (16)$$

The constraint "$I_m(\phi) = c_m$" is omitted when $m = 0$ because of Eq.(4). The pair (13)-(14) gives upper and lower bounds of the cumulant subject to the uncertainty size $\varepsilon$ evaluated in terms of the divergence $\mathbb{D}(\phi)$. By definitions of $\overline{I}_{k,m}, \underline{I}_{k,m}$, we have

$$0 \leq \underline{I}_{k,m} \leq \mathrm{Cum}_k(1) \leq \overline{I}_{k,m}, \quad (17)$$

showing that $\overline{I}_{k,m}, \underline{I}_{k,m}$ sandwich the cumulant $\mathrm{Cum}_k(1)$ of the benchmark case. The goal of the optimization problems is to find the maximizing/minimizing $\phi = \phi^*$, which corresponds to the worst-case uncertainties in the supOU process given an uncertainty size.

***Remark 4*** Because the estimation of a cumulant would be more sensitive to data for higher $m$, as discussed above, it is reasonable to assume that $k > m$. In **Section 4**, we focus on the cases $(k,m)=(1,0)$ (the average is uncertain) and $(k,m)=(2,1)$ (the variance is uncertain, and the average is more reliable).

***Remark 5*** The relationship (17) shows a qualitative difference between $\overline{I}_{k,m}, \underline{I}_{k,m}$ because the latter is not bounded from above, whereas the former is bounded from below by 0 and from above by $\mathrm{Cum}_k(1)$. The **Upper-bound case** is therefore more delicate.

### 3.2 Mathematical analysis

We study the existence and computability of optimizing $\phi$ in Eqs.(13)-(14). In the sequel, we assume that $k > m \geq 1$ because the results for the case $m = 0$ follow by suitably omitting the constraint in Eq.(15) and **Remark 4**. We first analyze the **Upper-bound case**, which is more delicate, and then the **Lower-bound case**.

#### 3.2.1 Upper-bound case

We need to carefully formulate the divergence because otherwise $\overline{I}_{k,m} = +\infty$. We analyze the **Upper-bound case** by focusing on when this unboundedness issue occurs with the help of Musielak–Orlicz spaces. For this purpose, we set the following function that is convex and increasing for $\phi_0 \geq 0$ by **Lemma 1**:

$$\overline{\Phi}(w_0,\alpha_0,\phi_0) = \begin{cases} 0 & (0 \leq \phi_0 \leq 1) \\ \Phi(w_0,\alpha_0,\phi_0) & (\phi_0 > 1) \end{cases}, \quad w_0,\alpha_0 > 0. \tag{18}$$

We set the convex conjugate $\Psi$ of $\Phi$ with respect to the third argument: for each $w_0,\alpha_0 > 0$,

$$\Psi(w_0,\alpha_0,\psi_0) = \sup_{\phi_0 > 0}\left\{\psi_0\phi_0 - \Phi(w_0,\alpha_0,\phi_0)\right\}, \quad \psi_0 \in \mathbb{R}. \tag{19}$$

Similarly, the convex conjugate of $\overline{\Phi}$ is denoted as $\overline{\Psi}$, which is given by Eq.(19) with the replacement $\Phi \Rightarrow \overline{\Phi}$. We obtain **Lemma 2** by direct calculations in Eq.(19) (see also Proposition 3.11 in Fröhlich and Williamson [36]):

***Lemma 2***

*For each* $w_0,\alpha_0 > 0$ *and* $\psi_0 \in \mathbb{R}$,

$$\Psi(w_0,\alpha_0,\psi_0) = \frac{w_0}{\alpha_0}\left\{\max\left\{0, 1+(\alpha_0-1)\frac{\psi_0}{w_0}\right\}^{\frac{\alpha_0}{\alpha_0-1}} - 1\right\} \tag{20}$$

*if* $\alpha_0 \neq 1$ *and*

$$\Psi(w_0,\alpha_0,\psi_0) = w_0\left(\exp\left(\frac{\psi_0}{w_0}\right) - 1\right) \tag{21}$$

*if* $\alpha_0 = 1$ *, where the right-hand side of Eq.(20) is understood to be* $+\infty$ *if* $0 < \alpha_0 < 1$ *and* $1 + (\alpha_0 - 1)\psi_0 \le 0$ *, and similarly,*

$$\frac{\partial \Psi(w_0, \alpha_0, \psi_0)}{\partial \psi_0} = \begin{cases} \max\left\{0, 1 + (\alpha_0 - 1)\dfrac{\psi_0}{w_0}\right\}^{\frac{1}{\alpha_0 - 1}} & (\alpha_0 \ne 1) \\ \exp\left(\dfrac{\psi_0}{w_0}\right) & (\alpha_0 = 1) \end{cases}, \tag{22}$$

*where the right-hand side of Eq.(22) is set to* $+\infty$ *if* $0 < \alpha_0 < 1$ *and* $1 + (\alpha_0 - 1)\psi_0 \le 0$ *. This* $\Psi$ *is nondecreasing and convex as a function of* $\psi_0 \in \mathbb{R}$ *. Moreover,* $\bar{\Psi} = \Psi$ *for* $\psi_0 \ge 0$ *. Finally,*

$$-\frac{w_0}{\alpha_0} \le \Psi(w_0, \alpha_0, \psi_0) \le \bar{\Psi}(w_0, \alpha_0, \psi_0), \quad \psi_0 \in \mathbb{R}. \tag{23}$$

The product measure $\pi\nu$ is not necessarily a probability measure on $(0, +\infty)^2$, but we can construct a probability measure on the basis of the constraint $I_m(\phi) = c_m$; we can rewrite it as follows:

$$\int_{r=0}^{r=+\infty} \int_{z=0}^{z=+\infty} \phi(r, z) \underbrace{\frac{z^m}{c_m m r} \pi(\mathrm{d}r) \nu(\mathrm{d}z)}_{p_m(\mathrm{d}r\mathrm{d}z)} = 1 \tag{24}$$

with

$$p_m(\mathrm{d}r\mathrm{d}z) = \left( \underbrace{\int_{r=0}^{r=+\infty} \int_{z=0}^{z=+\infty} \frac{z^m}{mr} \pi(\mathrm{d}r) \nu(\mathrm{d}z)}_{c_m} \right)^{-1} \frac{z^m}{mr} \pi(\mathrm{d}r) \nu(\mathrm{d}z), \quad r, z > 0 \tag{25}$$

being a probability measure on $(0, +\infty)^2$ by Eq.(16). Then, $\phi$ is seen as a Radon–Nikodym derivative that satisfies $\mathbb{E}_m[\phi] = 1$ , where $\mathbb{E}_m$ represents the expectation based on $p_m$ . In the sequel, we symbolically write an integration with respect to $p_m$ as an expectation $\mathbb{E}_m$ in a suitable probability space; e.g., $\mathbb{E}_m[\phi]$ can be identified as the left-hand side of Eq.(24).

With $p_m$, we can reformulate the optimization problem (13) as follows:

$$\text{Find } \bar{I}_{k,m} = \frac{m c_m}{k} \sup_{\phi \in \mathfrak{A}_{m,\varepsilon}} \mathbb{E}_m\left[z^{k-m} \phi(r, z)\right] \tag{26}$$

with $\mathfrak{A}_{m,\varepsilon}$ reformulated as a

$$\mathfrak{A}_{m,\varepsilon} = \left\{ \phi : (0, +\infty)^2 \to [0, +\infty) \,\middle|\, \mathbb{E}_m[\phi] = 1,\ \mathbb{D}'(\phi) \le \frac{\varepsilon}{m c_m} \right\}, \tag{27}$$

where

$$\mathbb{D}'(\phi) = \int_{r=0}^{r=+\infty} \int_{z=0}^{z=+\infty} \frac{r}{z^m} \Phi\big(w(r,z), \alpha(r,z), \phi(r,z)\big)\, p_m(\mathrm{d}r\mathrm{d}z) = \mathbb{E}_m\left[\Phi(w', \alpha, \phi)\right]. \tag{28}$$

We omit the arguments of $\Phi$ in Eq.(28). In Eq.(28), we set $w' = rz^{-m}w$, which is an inhomogeneous and hence state-dependent function even when $w$ is a constant. Now, the optimization problem has been reformulated as a maximization problem to determine a maximizing Radon–Nikodym derivative. The price to pay is that the weight $w'$ may be more complex than $w$.

By **Lemma 1** and the reformulation above, we can associate a Musielak–Orlicz space, which is a Banach space (Section 2 in Strati [32]), with the function $\bar{\Phi}$.

***Definition 1***

*The Musielak–Orlicz space* $L_{\bar{\Phi}}$ *is a collection of all real random variables* $\phi$ *on* $(0,+\infty)^2$ *such that*

$$\mathbb{E}_m\left[\bar{\Phi}\left(w',\alpha,\frac{|\phi|}{c}\right)\right] < +\infty \tag{29}$$

*for some constant* $c>0$, *and the space* $L_{\bar{\Phi}}$ *is equipped with the norm*

$$\|\phi\|_{\bar{\Phi}} = \inf\left\{c>0 \middle| \mathbb{E}_m\left[\bar{\Phi}\left(w',\alpha,\frac{|\phi|}{c}\right)\right] \le 1\right\}. \tag{30}$$

We associate another Musielak–Orlicz space $L_{\bar{\Psi}}$ with $\bar{\Psi}$, which is defined with the formal replacement $\bar{\Phi} \Rightarrow \bar{\Psi}$ in **Definition 1**. The Musielak–Orlicz space $L_{\bar{\Phi}}$ is reduced to a classical Orlicz space if $w$ and $\alpha$ are constant functions.

We consider the following standard structural condition about $\bar{\Phi}$, the so-called (uniform) $\Delta_2$ condition, which ensures the reflexivity of $L_{\bar{\Phi}}$ with its dual $L_{\bar{\Psi}}$ (e.g., Section 2 in Strat [32]).

***Assumption 3*** *There exists a constant* $K>0$ *such that*

$$\bar{\Phi}(w_0,\alpha_0,2\phi_0) \le K\bar{\Phi}(w_0,\alpha_0,\phi_0) \ \textit{for all} \ w_0,\alpha_0>0 \ \text{and} \ \phi_0>K. \tag{31}$$

The next proposition gives an existence result and a dual reformulation of the **Upper-bound case**, which is computationally accessible. Musielak–Orlicz spaces play a role here.

***Proposition 3***

*Assume that* $z^{k-m} \in L_{\bar{\Psi}}$ *and **Assumption 3**. The optimization problem* (13) *is reformulated as follows:*

$$\bar{I}_{k,m} = \frac{mc_m}{k}\inf_{\mu\in\mathbb{R},\ \tau>0}\left\{\frac{\varepsilon}{mc_m}\tau+\mu+\tau\mathbb{E}_m\left[\Psi\left(w'(r,z),\alpha(r,z),\frac{z^{k-m}-\mu}{\tau}\right)\right]\right\}\left(=\frac{mc_m}{k}\inf_{\mu\in\mathbb{R},\ \tau>0}F(\tau,\mu)\right) \tag{32}$$

*with the infimum being attained at* $(\hat{\tau},\hat{\mu}) \in (0,+\infty)\times[0,+\infty)$. *Moreover, an optimal* $\phi=\phi^*$ *is given by*

$$\phi^*(r,z) = \frac{\partial \Psi}{\partial \psi}\left( w'(r,z), \alpha(r,z), \frac{z^{k-m} - \hat{\mu}}{\hat{\tau}} \right), \quad (33)$$

*and satisfies the condition **(I$_k$)**. The couple* $(\hat{\tau}, \hat{\mu})$ *solves the equations*

$$\frac{\partial F(\tau,\mu)}{\partial \tau} = 0 \quad \textit{and} \quad \frac{\partial F(\tau,\mu)}{\partial \mu} = 0 . \quad (34)$$

The partial derivative $\frac{\partial \Psi}{\partial \psi}$ is taken with respect to the third argument. Here, the $\Delta_2$ condition (**Assumption 3**) ensures the existence of a minimizer of Eq.(32). This condition is satisfied for Eq.(12) with strictly bounded $\alpha$, but the condition $z^{k-m} \in L_{\overline{\Psi}}$ may fail if $\alpha(r,z) < 1$ for some $(r,z)$ because of Eq.(20) unless $w$ is chosen properly (see **Section 3.2.4**). Note also that the assumption $z^{k-m} \in L_{\overline{\Psi}}$ can be checked once a divergence is given before the optimization problem is solved, while the condition **(I$_k$)** cannot be checked because it uses the optimizer $\phi^*$, which is unknown. However, if $z^{k-m} \in L_{\overline{\Psi}}$, then the condition **(I$_k$)** with $\phi^*$ holds true as stated in the proposition, and this is the place where Musielak–Orlicz spaces determine the existence of the cumulants of the supOU process (1).

***Remark 6*** Inspecting **Proof of Proposition 3** shows that this proposition holds true if $\Phi$ is replaced by another function that satisfies certain properties (nonnegativity, convexity, minimized at $\phi_0 = 1$, and $\Phi(\cdot,\cdot,+\infty) = +\infty$) in **Lemma 1**, the $\Delta_2$ condition (**Assumption 3**), and the continuous differentiability of $\Psi(\cdot,\cdot,\psi_0)$ for all $\psi_0 \in \mathbb{R}$. One can apply the arguments in Proof of Proposition 5.5 and Proof of Theorem 6.5 in Dommel and Pichler [52].

### 3.2.2 Lower bound case

We have

$$\underline{I}_{k,m} = \inf_{\phi \in \mathfrak{A}_{m,\varepsilon}} I_k(\phi) = \inf_{\phi \in \mathfrak{A}_{m,\varepsilon}} \{ -I_k(-\phi) \} = -\sup_{\phi \in \mathfrak{A}_{m,\varepsilon}} I_k(-\phi) . \quad (35)$$

This is essentially a maximization problem analogous to the **Upper-bound case**, and we have the following proposition.

***Proposition 4***

*Assume that* $z^{k-m} \in L_{\overline{\Psi}}$ *and **Assumption 3**. The optimization problem (14) is reformulated as follows:*

$$\underline{I}_{k,m} = -\frac{mc_m}{k} \inf_{\mu \in \mathbb{R},\, \tau > 0} \left\{ \frac{\varepsilon}{mc_m}\tau + \mu + \tau \mathbb{E}_m \left[ \Psi\left( w'(r,z), \alpha(r,z), \frac{-z^{k-m} - \mu}{\tau} \right) \right] \right\} \left( = -\frac{mc_m}{k} \inf_{\mu \in \mathbb{R},\, \tau > 0} G(\tau,\mu) \right) \quad (36)$$

*with the infimum attained at* $(\hat{\tau}, \hat{\mu}) \in (0,+\infty) \times (-\infty, 0]$. *Moreover, an optimal* $\phi = \phi^*$ *is given by*

$$\phi^*(r,z) = \frac{\partial \Psi}{\partial \psi}\left( w'(r,z), \alpha(r,z), \frac{-z^{k-m} - \hat{\mu}}{\hat{\tau}} \right), \quad (37)$$

*and satisfies the condition **(I$_k$)**. The couple* $(\hat{\tau}, \hat{\mu})$ *solves the equations*

$$\frac{\partial G(\tau,\mu)}{\partial \tau} = 0 \quad \text{and} \quad \frac{\partial G(\tau,\mu)}{\partial \mu} = 0. \quad (38)$$

### 3.2.3 Ordering Musielak–Orlicz spaces

There is an ordering property among Musielak–Orlicz spaces, as shown in the next proposition, which is a state-dependent extension of the result for classical Orlicz spaces (Proposition 4 in Yoshioka and Yoshioka [30]). Here, the Musielak–Orlicz space associated with $\bar{\Phi}(w',\alpha,\cdot)$ is denoted as $L_{\bar{\Phi}(w',\alpha,\cdot)}$ with the associated norm $\|\cdot\|_{\bar{\Phi}(w',\alpha,\cdot)}$.

***Proposition 5***

*If* $\bar{\Phi}(w_1',\alpha_1,\cdot) \leq \bar{\Phi}(w_2',\alpha_2,\cdot)$, *then* $L_{\bar{\Phi}(w_2',\alpha_2,\cdot)} \subset L_{\bar{\Phi}(w_1',\alpha_1,\cdot)}$. *Moreover, if* $\phi \in L_{\bar{\Phi}(w_2',\alpha_2,\cdot)}$, *then* $\|\phi\|_{(w_2',\alpha_2,\cdot)} \geq \|\phi\|_{(w_1',\alpha_1,\cdot)}$.

The statement of this proposition is simple but intuitive. In the context of our optimization problems, a smaller Musielak–Orlicz space is more pessimistic about uncertainties because it corresponds to a greater divergence.

### 3.2.4 Specific setting

We discuss the implications of **Propositions 3-4,** focusing on a specific setting that is encountered in applied problems. Here, we use the measures $\pi$ and $\nu$ rather than $p_m$ because the latter has been introduced for technical reasons. Similarly, we use $w$ rather than $w'$. We study when the condition **(I$_k$)** is satisfied with the optimal $\phi = \phi^*$ and how the autocorrelation function would behave. We focus mainly on the **Upper bound case** in this subsection.

Because we are interested in long-memory phenomena in supOU processes, we assume the following gamma distribution for the reversion measure $\pi$, which is a nominal model for reversion speeds [e.g., 8,9,30]:

$$\pi(\mathrm{d}r) = c_\pi r^{A-1} e^{-\frac{r}{B}} \mathrm{d}r, \quad r > 0 \quad (39)$$

with the shape parameter $A > 1$, scale parameter $B > 0$, and normalization constant $c_\pi > 0$. For the Lévy measure $\nu$, we assume the following tempered-stable model, which is a common model for describing jumps [e.g., 20,53,54]:

$$\nu(\mathrm{d}z) = c_\nu \frac{e^{-pz}}{z^{q+1}} \mathrm{d}z, \quad z > 0 \tag{40}$$

with the shape parameter $q < 1$, tilting parameter $p > 0$, and intensity parameter $c_\nu > 0$. Under the present setting, the condition **(I*k*)** becomes

$$I_k = \frac{c_\pi c_\nu}{k} \int_{r=0}^{r=+\infty} \int_{z=0}^{z=+\infty} \phi(r,z) r^{A-2} z^{k-q-1} e^{-\frac{r}{B}-pz} \mathrm{d}r < +\infty, \quad k \in \mathbb{N}. \tag{41}$$

This implies that **(I*m*)** is satisfied if **(I*k*)** is ( $m < k$ ). We assume that $m < k$ in the rest of this subsection.

Now, we investigate the link between the singularity of the measures $\pi, \nu$ and the function $\Phi$. First, we assume that $\alpha(\cdot,\cdot) > 1$. By direct calculations,

$$\phi^*(r,z) = \max\left\{0, 1 + \frac{\alpha(r,z)-1}{w'(r,z)} \frac{z^{k-m} - \hat{\mu}}{\hat{\tau}}\right\}^{\frac{1}{\alpha(r,z)-1}} = \max\left\{0, 1 + \frac{\alpha(r,z)-1}{w(r,z)r} \frac{z^k - \hat{\mu} z^m}{\hat{\tau}}\right\}^{\frac{1}{\alpha(r,z)-1}}. \tag{42}$$

By (41), the admissibility of $\phi^*$ requires

$$\int_{r=0}^{r=+\infty} \int_{z=0}^{z=+\infty} \left\{\frac{\alpha(r,z)-1}{w(r,z)}\right\}^{\frac{1}{\alpha(r,z)-1}} r^{A-\frac{1}{\alpha(r,z)-1}-2} z^{\frac{\alpha(r,z)}{\alpha(r,z)-1}k-q-1} e^{-\frac{r}{B}-pz} \mathrm{d}r < +\infty. \tag{43}$$

If we additionally assume that $w$ is bounded and strictly positive, then Eq.(43) is reduced to

$$\int_{r=0}^{r=+\infty} \int_{z=0}^{z=+\infty} r^{A-\frac{1}{\alpha(r,z)-1}-2} z^{\frac{\alpha(r,z)}{\alpha(r,z)-1}k-q-1} e^{-\frac{r}{B}-pz} \mathrm{d}r < +\infty, \tag{44}$$

and the condition **(I*k*)** is satisfied if

$$A - \frac{1}{\alpha(r,z)-1} - 2, \frac{\alpha(r,z)}{\alpha(r,z)-1} k - q - 1 \geq -1 + \delta \tag{45}$$

for a sufficiently small constant $\delta > 0$. Because Eq.(45) is rewritten as follows (the second inequality can be trivially satisfied by choosing a suitably smaller $\delta$ when necessary):

$$\alpha(r,z) \geq 1 + \frac{1}{A-\delta-1} \quad \text{for} \quad 0 < \delta < A-1, \tag{46}$$

the generalized alpha divergence is suited to modeling uncertainties in the supOU process if $\alpha$ is sufficiently large. Because using a larger $\alpha$ results in a divergence that is less sensitive to uncertainties [30], one may use a state-dependent $\alpha$ by focusing on specific ranges of $(r,z)$. For example, in the context of modeling streamflow discharge, the focus may be high (large $z$ ) and/or prolonged (small $r$ ) flooding events. Similarly, one can also adaptively choose $w(r,z)$ under the same idea for each given function $\alpha$. If we assume that $\alpha$ is independent of $z$, then the argument analogous to that in **Proof of Proposition 2** predicts the tail behavior of the autocorrelation function as follows:

$$\mathrm{ACF}(h) = O\left(h^{-\left(A-\frac{1}{\alpha_{r=0}-1}-1\right)}\right) \quad \textit{for} \quad h >> 1, \tag{47}$$

showing that the ACF is predicted to be heavier in this case, and the heaviness becomes stronger for smaller $\alpha$ and hence for divergence closer to the KL one.

Second, we assume that $0<\alpha(\cdot,\cdot)<1$. This is possibly an ill-posed case where the condition **($\mathbf{I}_k$)** is never satisfied. Indeed, the admissibility of $\phi^*$ requires

$$\int_{r=0}^{r=+\infty}\int_{z=0}^{z=+\infty} \max\left\{0, 1-\frac{1-\alpha(r,z)}{w(r,z)r}\frac{z^k-\hat{\mu}z^m}{\hat{\tau}}\right\}^{\frac{-1}{1-\alpha(r,z)}} r^{A-2}z^{k-q-1}e^{-\frac{r}{B}-pz}\mathrm{d}r<+\infty \tag{48}$$

along with

$$1-\frac{1-\alpha(r,z)}{w(r,z)r}\frac{z^k-\hat{\mu}z^m}{\hat{\tau}}>0 \quad \text{for all} \quad r,z>0 . \tag{49}$$

The positivity condition Eq.(49) fails if the coefficient $\dfrac{1-\alpha(r,z)}{w(r,z)}$ does not vanish for large $z$ and small $r$. One may choose $w(r,z)$ so that this singularity is resolved while $\alpha(r,z)$ is strictly bounded between 0 and 1; however, this case needs a more complicated (and possibly difficult to interpret and extremely insensitive) modeling strategy of $w$ and $\alpha$, and therefore, our application in **Section 4** assumes $\alpha \geq 1$.

Finally, we move to the case $\alpha \equiv 1$, which corresponds to evaluating uncertainties by using a state-dependent KL divergence. By direct calculations, we have

$$\phi^*(r,z)=\exp\left\{\frac{1}{w(r,z)r}\frac{z^k-\hat{\mu}z^m}{\hat{\tau}}\right\}, \tag{50}$$

and the admissibility of $\phi^*$ requires

$$\int_{r=0}^{r=+\infty}\int_{z=0}^{z=+\infty}\exp\left\{\frac{1}{\hat{\tau}w(r,z)r}z^k\right\}r^{A-2}z^{k-q-1}e^{-\frac{r}{B}-pz}\mathrm{d}r<+\infty . \tag{51}$$

This condition is insightful because it always fails if we assume a classical KL divergence ( $w$ is a constant function) irrespective of $k$ because of the singular function $\dfrac{1}{r}$ in the exponential. This singularity issue for the classical KL divergence was reported previously [30]. We can resolve this issue if $k=1$ by choosing $w$ so that the function $\dfrac{1}{w(r,z)r}$ is bounded from above. The Musielak–Orlicz approach allows for such a flexible modeling strategy; e.g., one may set

$$\frac{1}{w(r,z)r}=\frac{1-e^{-Wr}}{cr} \quad \text{or equivalently} \quad w(r,z)=\frac{c}{1-e^{-Wr}} \tag{52}$$

with constants $c,W>0$. This $w$ gives a strictly bounded $1/(wr)\leq W$ and resolves the singularity issue caused by $1/r$ using the regularization near $r=0$. The strength of this regularization can be controlled by $W>0$ such that a larger $W$ indicates a weaker regularization effect and that there is no regularization under the limit $W\to+\infty$. One can also consider a $z$-dependent $w$ if it remains bounded and does not encounter the singularity issue at $r=0$, implying the versatility of the present approach based on Musielak–Orlicz spaces. In this case, because there is no singularity of $\phi^*$ at $r=0$, the tail behavior in **Proposition 2** does not change, which is the case for our computation in **Section 4**. We computationally find that $\hat{\tau}$ is sufficiently large such that Eq.(51) is satisfied.

We briefly discuss the **Lower bound case** where the KL divergence does not fail because of

$$\phi^*(r,z)=\exp\left\{\frac{-1}{w(r,z)r}\frac{z^k-\hat{\mu}z^m}{\hat{\tau}}\right\}, \tag{53}$$

which satisfies the condition **(I$_k$)**; hence, we can choose $w$ as a constant in this case. The tail of the autocorrelation function of the distorted model does not change with $w$ given by Eq.(52). If $\alpha(\cdot,\cdot)>1$,

$$\phi^*(r,z)=\max\left\{0,-\frac{\alpha(r,z)-1}{w(r,z)r}\frac{z^k-\hat{\mu}z^m}{\hat{\tau}}+1\right\}^{\frac{1}{\alpha(r,z)-1}}, \tag{54}$$

and the condition **(I$_k$)** is always satisfied if $z^{k-m}\in L_\Psi$ because $\hat{\mu}\le 0$.

### 3.3 Computational method

By **Propositions 4-5**, we can solve the optimization problems by using a gradient descent method for the Lagrangian multipliers $\mu,\tau$. We aim to present a simple computational method for our optimization problems but do not compare algorithms.

For the **Upper-bound case**, we want to find a solution $(\hat{\tau},\hat{\mu})$, optimal Lagrangian multipliers solving Eq.(34), where by elementary calculations,

$$\frac{\partial F(\tau,\mu)}{\partial\tau}=\frac{\varepsilon}{mc_m}-\mathbb{E}_m\left[\Phi\left(w',\alpha,\frac{\partial\Psi}{\partial\psi}\left(w',\alpha,\frac{z^{k-m}-\mu}{\tau}\right)\right)\right], \tag{55}$$

$$\frac{\partial F(\tau,\mu)}{\partial\mu}=\tau\left\{1-\mathbb{E}_m\left[\frac{\partial\Psi}{\partial\psi}\left(w',\alpha,\frac{z^{k-m}-\mu}{\tau}\right)\right]\right\}. \tag{56}$$

A difficulty in solving the system (34) is that the partial derivatives $\frac{\partial F}{\partial\tau}$ and $\frac{\partial F}{\partial\mu}$ are not globally Lipschitz continuous and are often assumed to safely implement gradient descent methods and their variants [55-57]. The loss of the Lipschitz continuity of partial derivatives is due to the appearance of the term $\frac{z^{k-m}-\mu}{\tau}$. Nevertheless, because Eq.(34) admits an interior solution (**Proposition 3**), we have $(\hat{\tau},\hat{\mu})\in\left[L^{-1},L\right]\times\left[L^{-1},L\right]$ with some sufficiently large constant $L>0$, and the partial derivatives of $F$ are Lipschitz continuous in this domain. We implement a gradient method with momentum to solve Eq.(34), as shown in **Section A2 in Appendix,** with a convergence study.

We need to discretize the expectations appearing in each partial differential in Eq.(34) to fully implement the gradient descent method. We apply the quantization method (Proof of Proposition 2 in Yoshioka [58]). Note that the probability measure $p_m$ is decomposable to the two probability measures $\pi_m = c_1 r^{-1}\pi$ and $\nu_m = c_2 z^m \nu$ with normalization constants $c_1, c_2$. Fix the computational resolution $M \in \mathbb{N}$, and set the quantiles $P_{i,M}$ and $Q_{i,M}$ that uniformly discretize $\pi_m, \nu_m$ as follows:

$$\frac{i}{M} = \int_0^{P_{i,M}} \pi_m(\mathrm{d}r) \quad \text{and} \quad \frac{i}{M} = \int_0^{Q_{i,M}} \nu_m(\mathrm{d}v), \quad i = 0,1,2,\ldots,M . \tag{57}$$

Afterward, we set $x_{i,M} = P_{2i-1,2M}$ and $y_{i,M} = Q_{2i-1,2M}$ ($i = 1,2,3,\ldots,M$). This quintile approximation has a pointwise error $M^{-1}$ between the discretized and true cumulative distributions of a probability measure that admits a density (Proof of Proposition 2 in Yoshioka [58]). We apply the following rule to discretize each expectation appearing in Eqs.(55)-(56): for generic $f : (0,+\infty)^2 \times \mathbb{R} \to \mathbb{R}$,

$$\mathbb{E}_m\left[ f\left( w'(r,z), \alpha(r,z), \frac{z^{k-m} - \mu}{\tau} \right) \right] \approx \frac{1}{M} \sum_{i,j=1}^{M} f\left( w'\left(x_{i,M}, y_{j,M}\right), \alpha\left(x_{i,M}, y_{j,M}\right), \frac{\left(y_{j,M}\right)^{k-m} - \mu}{\tau} \right). \tag{58}$$

Finally, we can numerically solve the **Lower bound case** in the same way by formal replacement $F \Rightarrow G$ in the explanation above, where

$$\frac{\partial G(\tau,\mu)}{\partial \tau} = \frac{\varepsilon}{mc_m} - \mathbb{E}_m\left[ \Phi\left( w', \alpha, \frac{\partial \Psi}{\partial \psi}\left( w', \alpha, \frac{-z^{k-m} - \mu}{\tau} \right) \right) \right], \tag{59}$$

$$\frac{\partial G(\tau,\mu)}{\partial \mu} = \tau \left\{ 1 - \mathbb{E}_m\left[ \frac{\partial \Psi}{\partial \psi}\left( w', \alpha, \frac{-z^{k-m} - \mu}{\tau} \right) \right] \right\}. \tag{60}$$

***Remark 7*** When $m = 0$, it suffices to formally fix $\mu = 0$ and simply find a minimizing $\tau$, where the functions $F$ and $G$ are convex. The mathematical and computational approaches for $m \geq 1$ carry over to this case with proper simplifications.

***Remark 8*** For large-scale problems, the expectations need to be estimated through some machine-learning scheme that is potentially more robust and efficient to better estimate them because the quantization method faces the curse-of-dimensionality.

## 4. Demonstrative application

### 4.1 Data description

The supOU process (1) and the optimization problems (13)-(14) are applied to real data. A novelty here is the modeling of uncertainties by using a state-dependent divergence of the KL type under a theoretical justification that was lacking in previous studies [30,31]. The computational results in **Section 4**, particularly those in **Section 4.3.2**, provide insights into how model uncertainties quantitatively affect flood and drought regimes and their durations at the study site under state-dependent uncertainties, connecting weighted divergences, Musielak–Orlicz processes, and long-memory processes.

We use the hourly streamflow discharge data of a Class-A River in a mountainous Hakusan area, Hakusan City, Ishikawa Prefecture, Japan (Kazarashi station: 36.16778 N, 136.63194 E, with an elevation of 494.0 (m)). Temporally fine (e.g., hourly to daily) streamflow discharge and its fluctuation have been suggested to exhibit long memory [e.g., 15,48]. River water in this area has served as a water resource for residents [59] and local industries such as Tofu production and the aquaculture of inland fish species. A dominant source of the long memory of the river water in the study area is considered to be shallow groundwater supply and snowpack accumulation, but detailed estimations of these sources have not been conducted thus far. Moreover, the study area has been suggested to experience an increased probability of annual flood occurrence under future climate change [60,61] and an increase in the temporal variability of streamflow discharge under regional warming scenarios [62], motivating us to investigate supOU processes subject to model uncertainties.

Hourly discharge data have been measured at Kazarashi station by the Ministry of Land, Infrastructure, Transport and Tourism; however, recently, most of the data are "missing", at least from June to July in 2023, probably because of technical reasons such as the failure of observation equipment[2], and future changes in streamflow characteristics under climate change are considered important. We collected hourly discharge data from April 1, 2016, to March 31, 2023, at Kazarashi station[2] in consideration of data availability (most parts of the data are "missing" in January to March in 2016, and in total 61,344 data points exist and 577 of them are missing). **Figure 1** shows the streamflow discharge data used in this paper. **Figure 2** compares the empirical and theoretical autocorrelation functions of the discharge data. We identified a classical supOU process with reversion and Lévy measures by using Eqs.(39)-(40) and a moment-matching method that has already been validated in previous studies [e.g., 63], along with a comparison between empirical and theoretical cumulants (**Table 1**). As shown in **Table 1**, the streamflow discharge at this station has a long-memory property because the parameter $A$, the Hurst index, is approximately 1.5, and this characteristic is captured by the theory (**Figure 2**).

**Table 1.** Theoretical parameter values cumulants. The values inside "(, )" represent empirical cumulants and their relative errors.

| $A$ (-) | 1.502.E+00 |
| --- | --- |
| $B$ (1/h) | 1.282.E-01 |
| $c_v$ ($\mathrm{m}^{3q}/\mathrm{s}^{q}/\mathrm{h}$) | 2.856.E-03 |
| $p$ ($\mathrm{s/m^3}$) | 7.474.E-03 |
| $q$ (-) | -2.429.E-01 |
| Average ($\mathrm{m^3/s}$) | 1.771.E+01 (1.757.E+01, 7.899.E-03) |
| Variance ($\mathrm{m^6/s^2}$) | 1.472.E+03 (1.478.E+03, 3.996.E-03) |
| Skewness (-) | 5.214.E+00 (4.977.E+00, 4.750.E-02) |
| Kurtosis (-) | 4.421.E+01 (4.518.E+01, 2.135.E-02) |

[2] Water Information system. http://www1.river.go.jp/cgi-bin/SiteInfo.exe?ID=304111284412020 (In Japanese. Last Accessed on February 24, 2026)

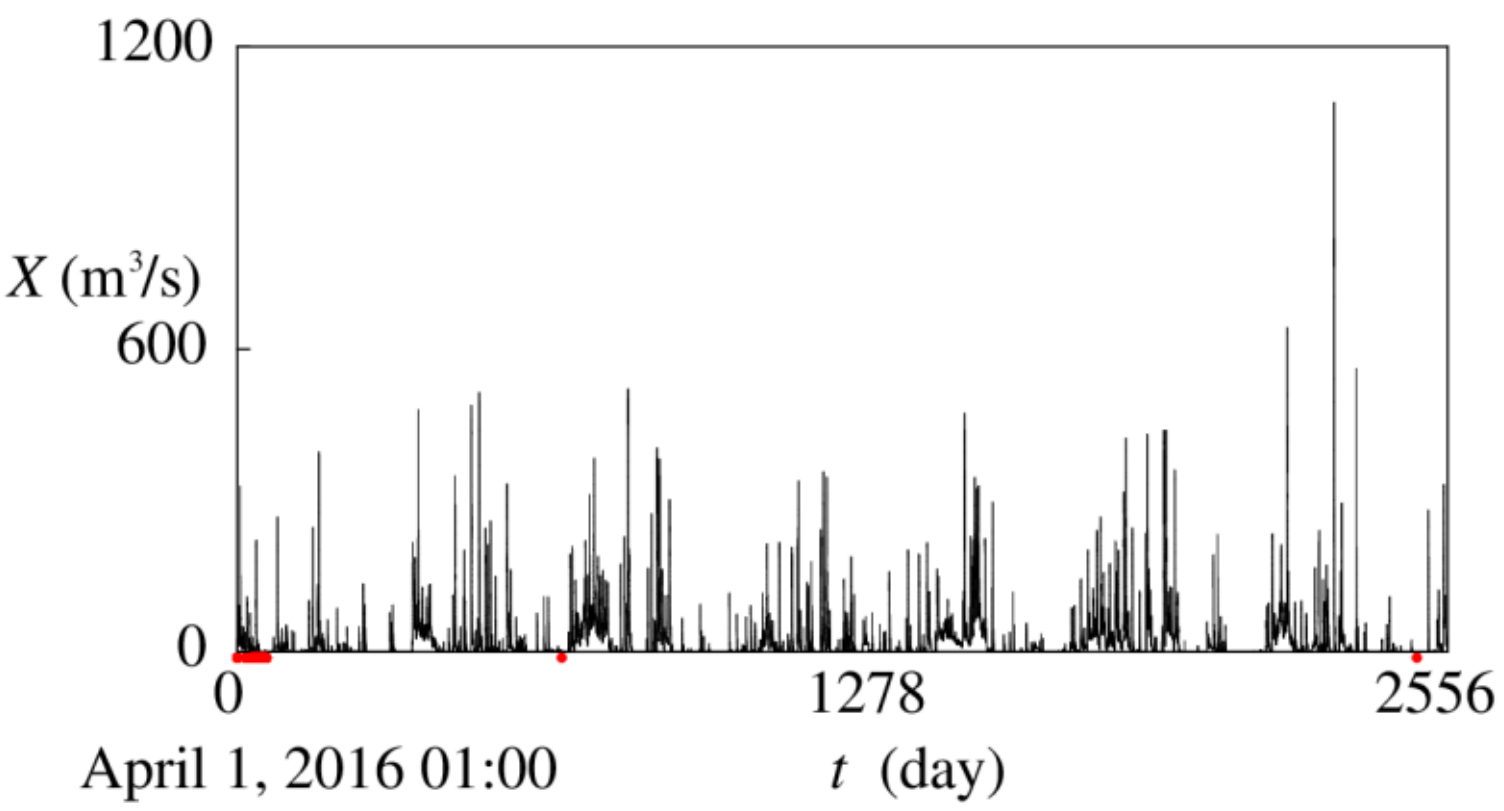

**Figure 1.** Discharge data at Kazarashi station. Red circles represent the time points with missing data.

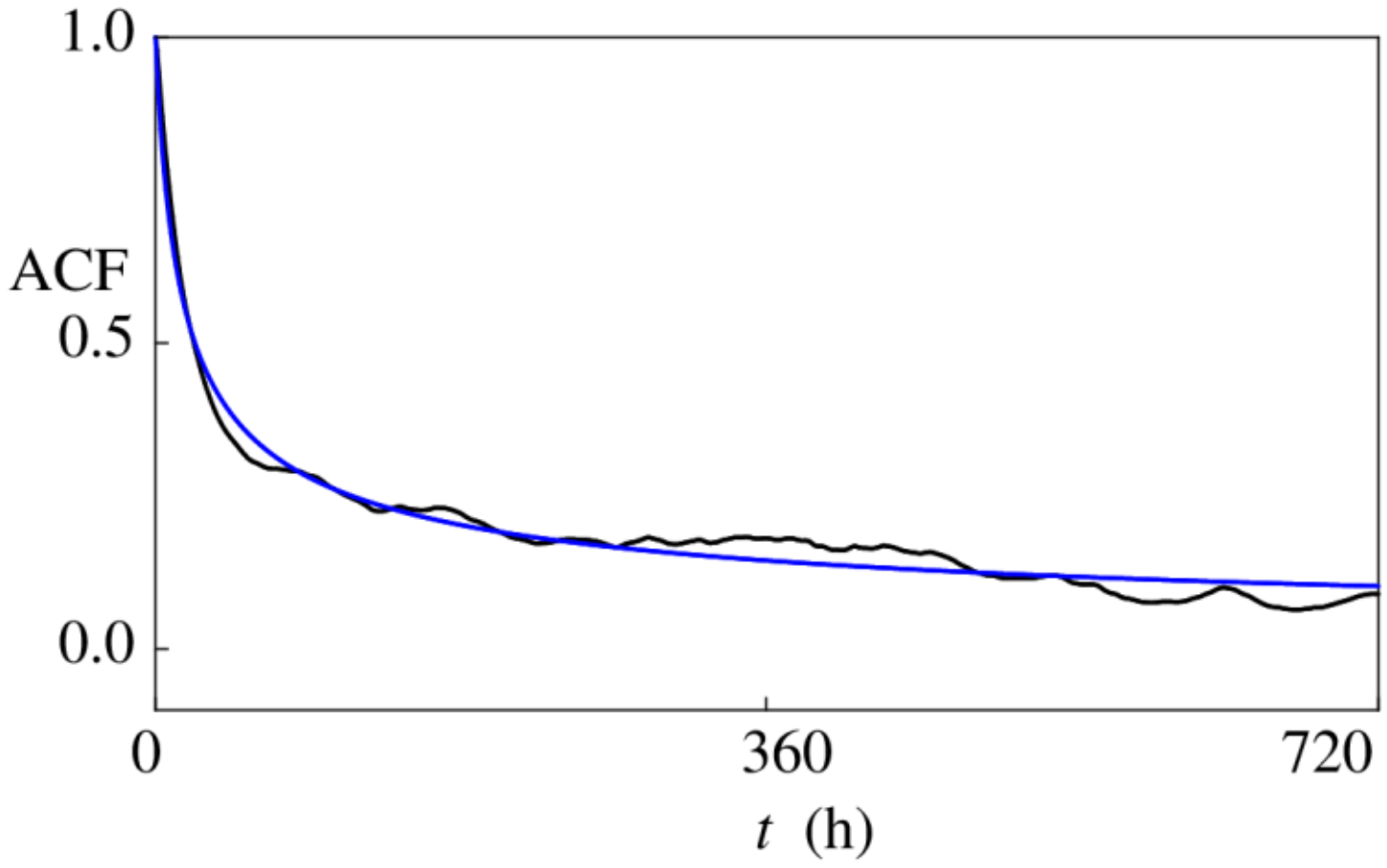

**Figure 2.** Autocorrelation function (ACF) at Kazarashi station. The black curve represents the empirical ACF, and the blue curve represents the theoretical ACF.

### 4.2 Computational setting

We compute the optimization problems in the **Upper- and Lower-bound cases** for different $\Phi$. We first use the following alpha divergence function:

$$\Phi\left(w(r),\alpha(r),x\right)=w(r)\frac{x^{\alpha(r)}-\alpha(r)(x-1)-1}{\alpha(r)\left(\alpha(r)-1\right)},\quad r,z>0\ \text{ and }\ x>0, \tag{61}$$

where $\alpha(\cdot)>1$ is a constant or continuous function of $r$ to account for the importance of reversion speeds, namely, timescales. For the visualization of $\Phi$ and its convex conjugate $\Psi$, see **Figures 3-4**. The weight $w$ is given by Eq.(52) so that the singularity at $r=0$ is regularized because some computational instability is faced in the gradient descent method, which is considered because of the appearance of very large expectants in Eqs.(55)-(56). The constant $c>0$ is determined to normalize the weight as $\int_0^{+\infty}\int_0^{+\infty}w'(r,z)\pi(\mathrm{d}r)\nu(\mathrm{d}z)=1$. The weight $w$ in Eq.(52) is valid if the memory structure of $X$ is not critically affected by uncertainties, which is assumed in **Section 4** for modeling simplicity. Even in this simplified setting, $w'(r,z)$ depends on both $r$ and $z$. We investigate the behavior of the Lagrangian multipliers and solutions to the optimization problems under (61). For this case, we set $(m,k)=(1,2)$ and examine $\varepsilon = 10^{2-0.6i}$ ($i = 0,1,2,...,10$), and the error threshold for convergence is $\delta = 10^{-5}$.

We also consider the following weighted KL divergence:

$$\Phi(w(r),x) = w(r)(x\ln x - x + 1),\quad r,z>0 \text{ and } x>0, \tag{62}$$

where $w$ is given by Eq.(52) and $c>0$ is the normalization constant. This case is interesting by itself, and moreover, it gives the following interpretation of the reversion and Lévy measures:

$$\phi^*(r,z)\pi(\mathrm{d}r)\nu(\mathrm{d}z) = \underbrace{\pi(\mathrm{d}r)}_{\substack{\text{Reversion measure}\\\text{(Not modified)}}} \times \underbrace{\exp\left(\pm\frac{c(1-e^{-Wr})}{\hat{\tau}r}z\right)\nu(\mathrm{d}z)}_{r\text{-dependent tempered-stable Lévy measure}},\quad r,z>0, \tag{63}$$

where the distorted Lévy measure is again of the tempered-stable type that is now inhomogeneous, which is computationally advantageous in our case; **Table 1** suggests that the Lévy measure $\nu$ in the present case is of the compound Poisson type with jumps distributed according to a gamma distribution. Therefore, with Eq.(63), jumps can be generated by gamma distributions depending on $r$. We compute the supOU process by using a discrete superposition scheme, which itself is not new but simple to implement. We discretize Eq.(1) for $r>0$ as follows on the basis of the representation formula of supOU processes (e.g., Proposition 1 in Fasen and Klüppelberg [51]; Eqs. (63)-(64) in Yoshioka [58]):

$$X_t \approx \sum_{i=1}^{M} X_t^{(i)} \text{ with } \mathrm{d}X_t^{(i)} = -r_i X_t^{(i)}\mathrm{d}t + \mathrm{d}Z_t^{(i)}, \tag{64}$$

where each $Z^{(i)}$ is a mutually independent jump process with Lévy measures $\frac{1}{B(A-1)} r_i M^{-1}\phi^*(r_i,z)\nu(\mathrm{d}z)$, with each $r_i$ being generated by the quantization in Eq.(57) with the formal replacement $A \to A-1$ in Eq.(39) and suitably update the normalization constant $c_\pi$, and each $X^{(i)}$ is computed using a common Euler–Maruyama time discretization, and jumps are generated by an acceptance-rejection method. For the present case, we set $(m,k)=(0,1)$ and examine $\varepsilon = 10^{1-0.6i}$ ($i = 0,1,2,...,10$), and the error threshold for convergence is $\delta = 10^{-7}$. We compute the threshold crossing dynamics of the supOU process at different $\varepsilon$ as explained later. The time increment for sampling each supOU process is 1/100 (h), and we simulate each sample path for 2,000 (yr) with a burn-in period of 2000 (yr) with a time increment of 1 (h).

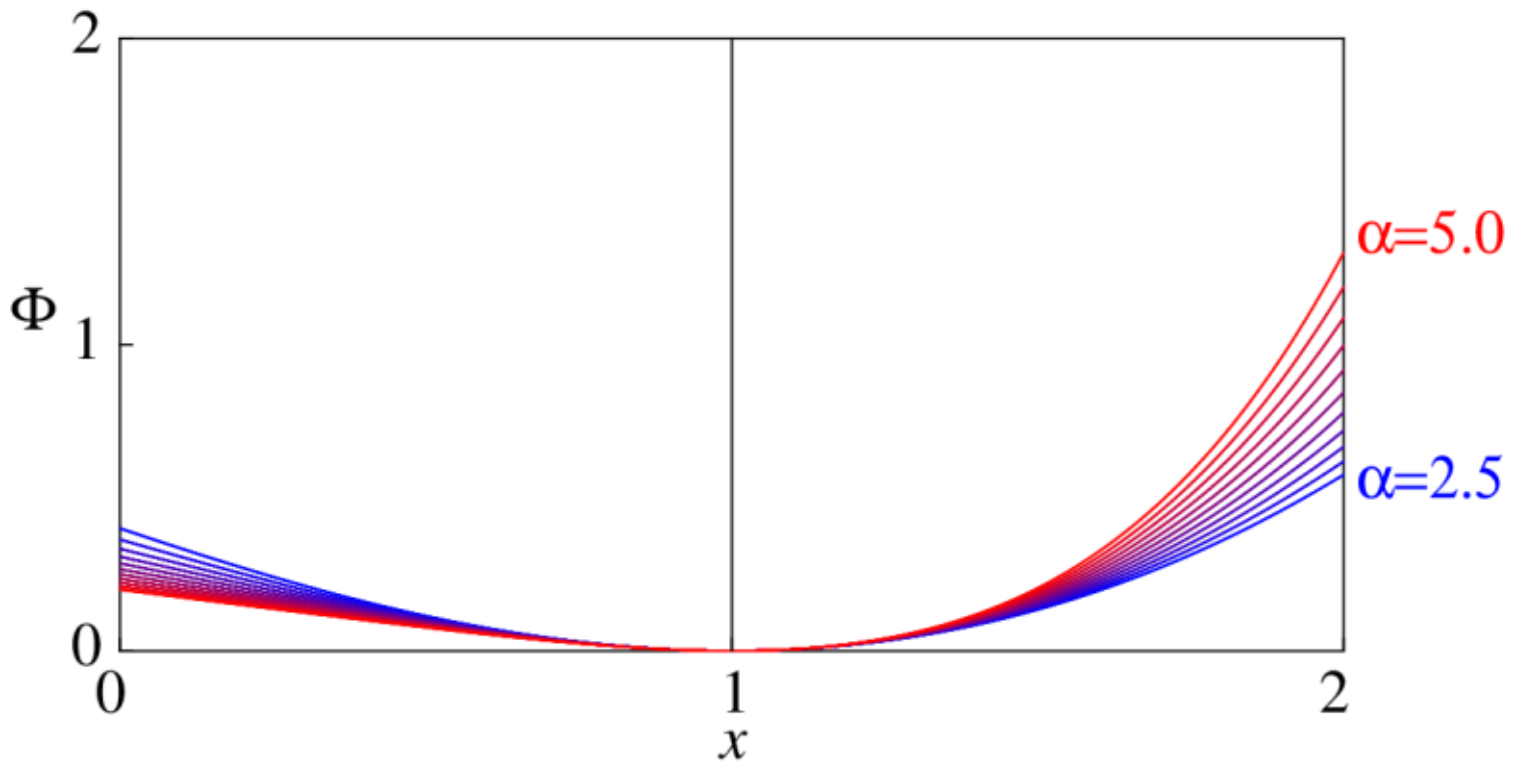

**Figure 3.** Homogeneous ( $w'$ is replaced by 1 and $\alpha$ by a constant) divergence function $\Phi = \Phi(x)$ for $\alpha = 2.5i$ ( $i = 0,1,2,...,10$ ).

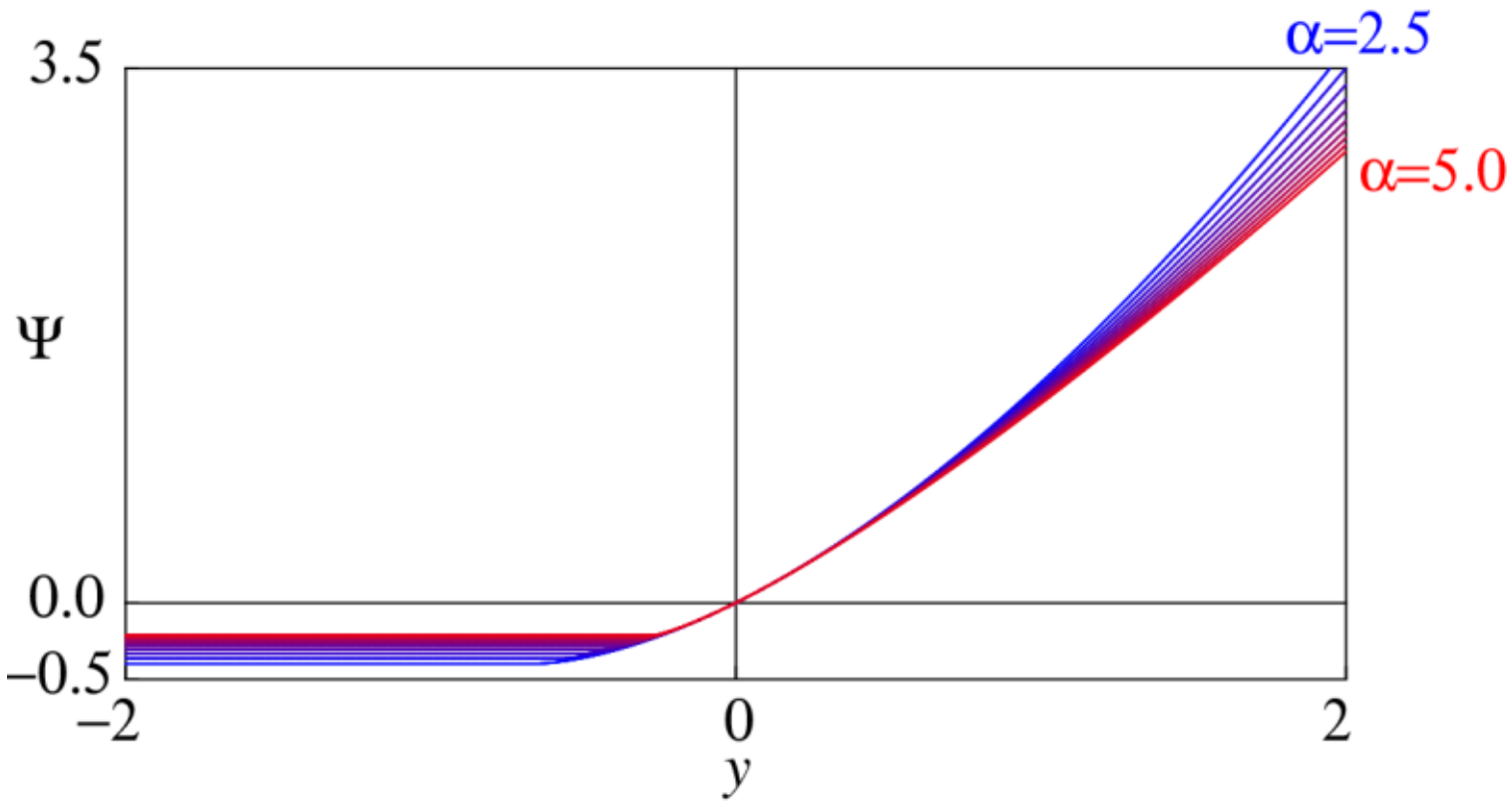

**Figure 4.** Convex divergence $\Psi = \Psi(y)$ of the homogeneous $\Phi$ in **Figure 3** for $\alpha = 2.5i$ ( $i = 0,1,2,...,10$ ).

### 4.3 Results and discussion

#### 4.3.1 Parameter dependence of the optimization problems

We examine three constant $\alpha$ cases $\alpha \equiv 2.5,\ 4.0,\ 5.0$ and two state-dependent $\alpha$ cases

$$\alpha(r) = 2.5 + (5.0 - 2.5)\min\{1, r\} \quad \text{and} \quad \alpha(r) = 5.0 + (2.5 - 5.0)\min\{1, r\}, \tag{65}$$

which are called the increasing and decreasing cases, respectively. We choose the regularization parameter value of $W = 10$ (h), and the influences of $W$ are discussed in **Section A3** in **Appendix**.

**Figures 5-7** show the computed variances, $\hat{\tau}$, and $\hat{\mu}$ for each case. **Figure 5** for the variance under distortion shows that, in the **Upper-bound case**, specifying smaller $\alpha$ values yields a more pessimistic overestimation of the variance at each $\varepsilon$. For the inhomogeneous $\alpha$, their results are between those of $\alpha \equiv 2.5,\ 5.0$, possibly because of Eq.(65), but the increasing and decreasing cases give distinctive results. The increasing case (resp., decreasing case) is more pessimistic for smaller (resp., larger) reversion speeds $r$ and gives results closer to those of $\alpha \equiv 2.5$ (resp., $\alpha \equiv 5.0$ ). A smaller $r$ corresponds to slowly decaying streamflow discharge, such as groundwater recharge to the river (e.g., baseflow), where the weight $w'$ becomes large. Evaluation of the state-dependence, which is based on Musielak–Orlicz spaces, allows for this inhomogeneity in $r$. The computational results for the **Lower bound case** are not significantly different among the cases examined here, suggesting higher sensitivity in the **Upper bound case** against $\alpha$ and hence the corresponding Musielak–Orlicz spaces.

**Figure 6** for the optimal Lagrangian multiplier $\hat{\tau}$ about the constraint of divergence shows that they are decreasing in $\varepsilon$ for all computational cases. This implies that perturbing the divergence constraint affects the objective more for larger $\varepsilon$ values, implying that tightening the uncertainty set rapidly decreases the objective for larger uncertainty sets (Chapter 5.6.2 in Boyd and Vandenberghe [44]). **Figure 7** for the optimal Lagrangian multiplier $\hat{\mu}$ about the constraint of average shows that $\hat{\mu}$ is not monotone with respect to $\varepsilon$; for the **Upper-bound case**, the monotonicity appears at approximately $\varepsilon$ values of 1 to 10, corresponding to the variance values doubling the benchmark value. For the **Lower-bound case**, the monotonicity appears at approximately $\varepsilon$ values of 0.001 to 0.1, with the variance values closer to (difference is smaller than a few percent of) the benchmark value. The reasons for these nonmonotone profiles are not theoretically clear but are considered due to the problem setting, except for the divergences because they commonly appear in the computational cases examined here. For $\hat{\mu}$, the data in **Figure 6** suggest that perturbing the constraint about the average discharge has greater effects on the objective for smaller $\varepsilon$ values, which is in contrast to the role of $\hat{\tau}$.

The conditional cumulative distributions of $p_m$ for $r$ and $z$ values are denoted as $Fr$ and $Fz$, which range from 0 to 1 as $r$ and $z$ increase from 0 to $+\infty$, respectively. **Figures 8-9** show the computed optimizers $\phi^*$ for the increasing case for the **Upper- and Lower-bound cases**, respectively. Decreasing $\varepsilon$ results in $\phi^*$ being closer to the constant function $\phi^* \equiv 1$ because a smaller distortion is allowed. For the **Upper bound case**, as shown in **Figure 8**, $\phi^*$ is large for large $Fz > 0.8$, which corresponds to the tail probability of jumps, implying that model uncertainties are judged because of the distortion of tails in the proposed mathematical framework. This implies that the pessimistic increase in streamflow discharge is due to the increase in the magnitude and frequency of large flood events rather than small flood events. Moreover, higher $\phi^*$ values for large $z$ and small $r$ values imply that the overestimation of baseflow events also contributes to the increase in streamflow variability; this point is consistent with the discussion for **Figure 5** of the variance under uncertainty. In contrast, for the **lower-bound case**, as shown in **Figure 9**, $\phi^*$ reaches large values for small $Fz$, suggesting that the pessimistic decrease in streamflow discharge is due to the decrease in base flows and flood events. In **Figures 8-9**, the heterogeneity of $\phi^*$ for the reversion speed $r$ seems rather weak compared with that of jumps $z$. For both the **Upper- and Lower-bound cases**, the optimal $\phi^*$ has more localized profiles for larger $\varepsilon$ values. The analysis of $\phi^*$ thus provides how the worst-case model uncertainties distort the reversion and Lévy measures.

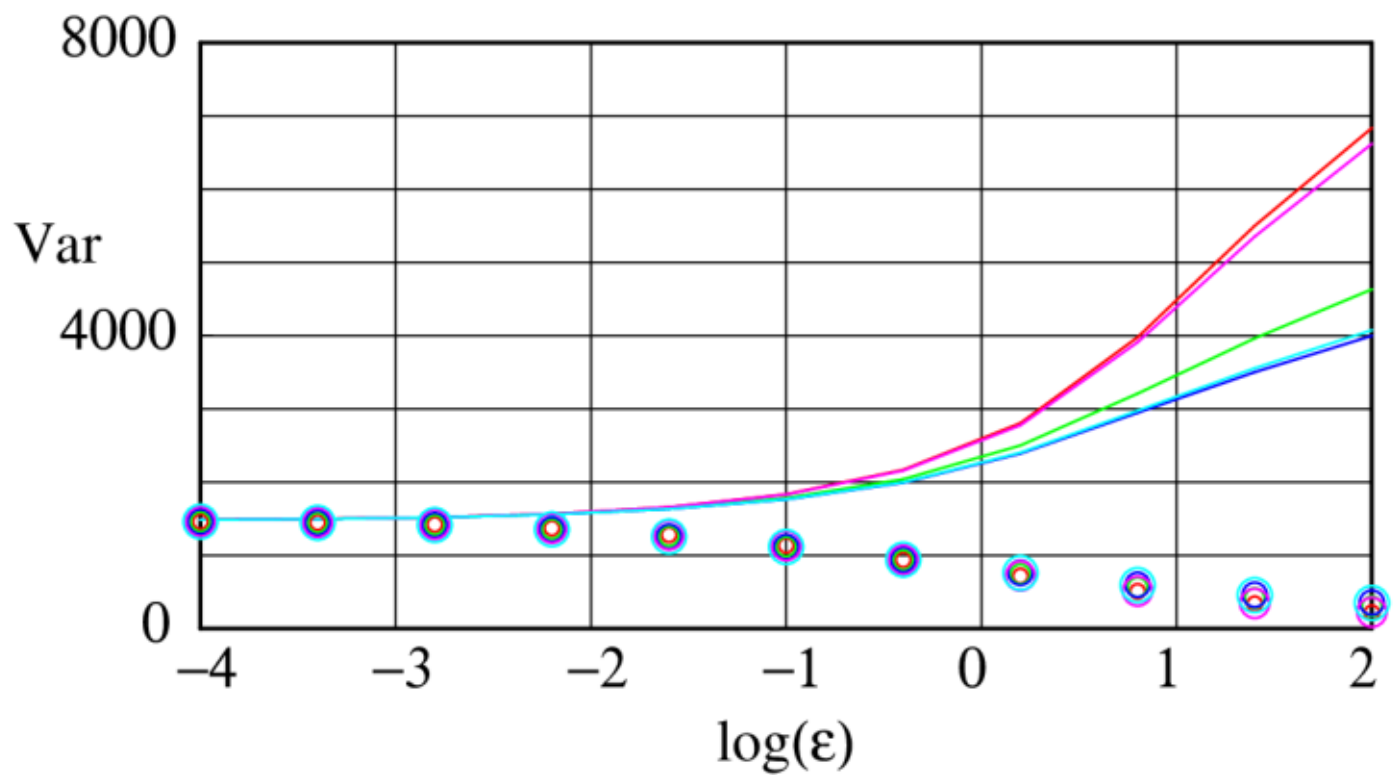

**Figure 5.** Computed variance (Var) ($m^6/s^2$) values of the distorted supOU processes for different values of $\varepsilon$ in **Upper-bound case** (curves) and **Lower-bound case** (circles). Colors represent $\alpha \equiv 2.5$ (red), $\alpha \equiv 4.0$ (green), $\alpha \equiv 5.0$ (blue), increasing $\alpha$ (magenta), and decreasing $\alpha$ (sky blue).

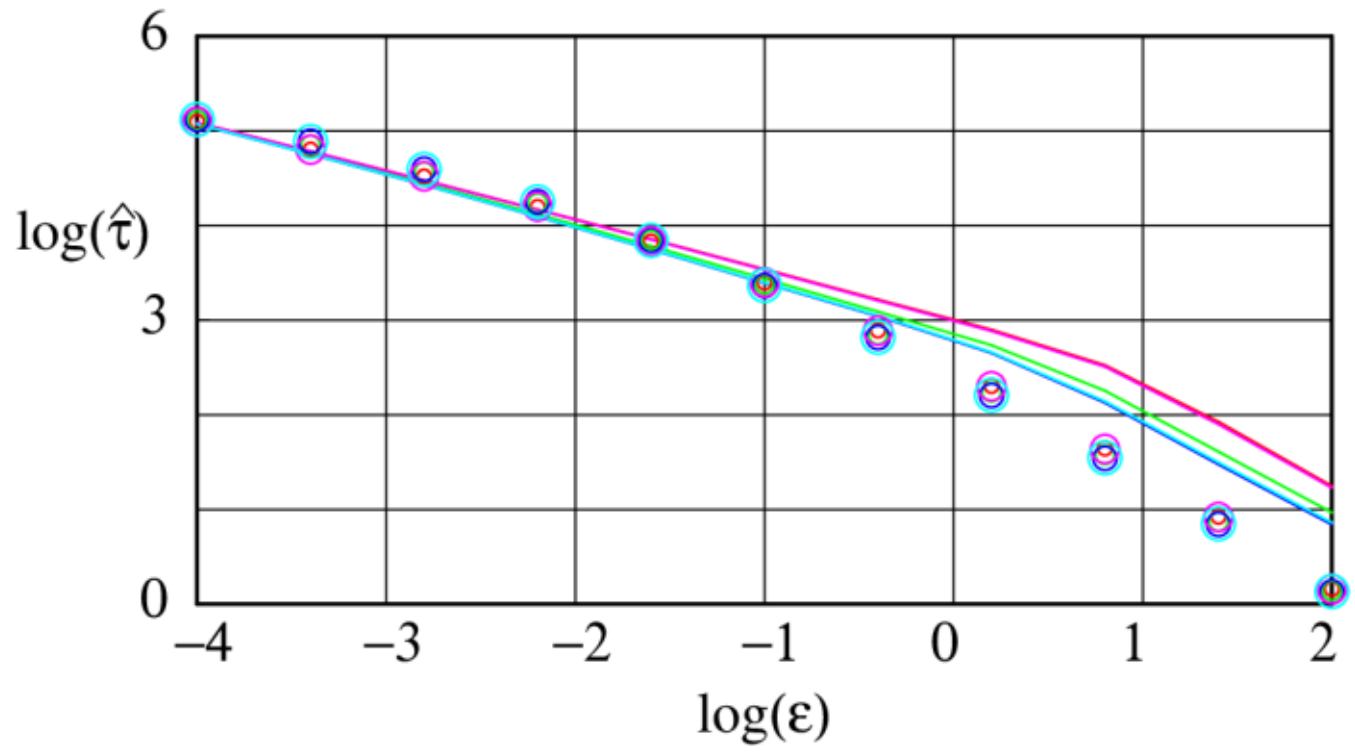

**Figure 6.** Computed $\tau = \hat{\tau}$ ($m^6/s^2$) values of the distorted supOU processes for different values of $\varepsilon$ in **Upper-bound** and **Lower-bound cases**. Same color legend as in **Figure 5**.

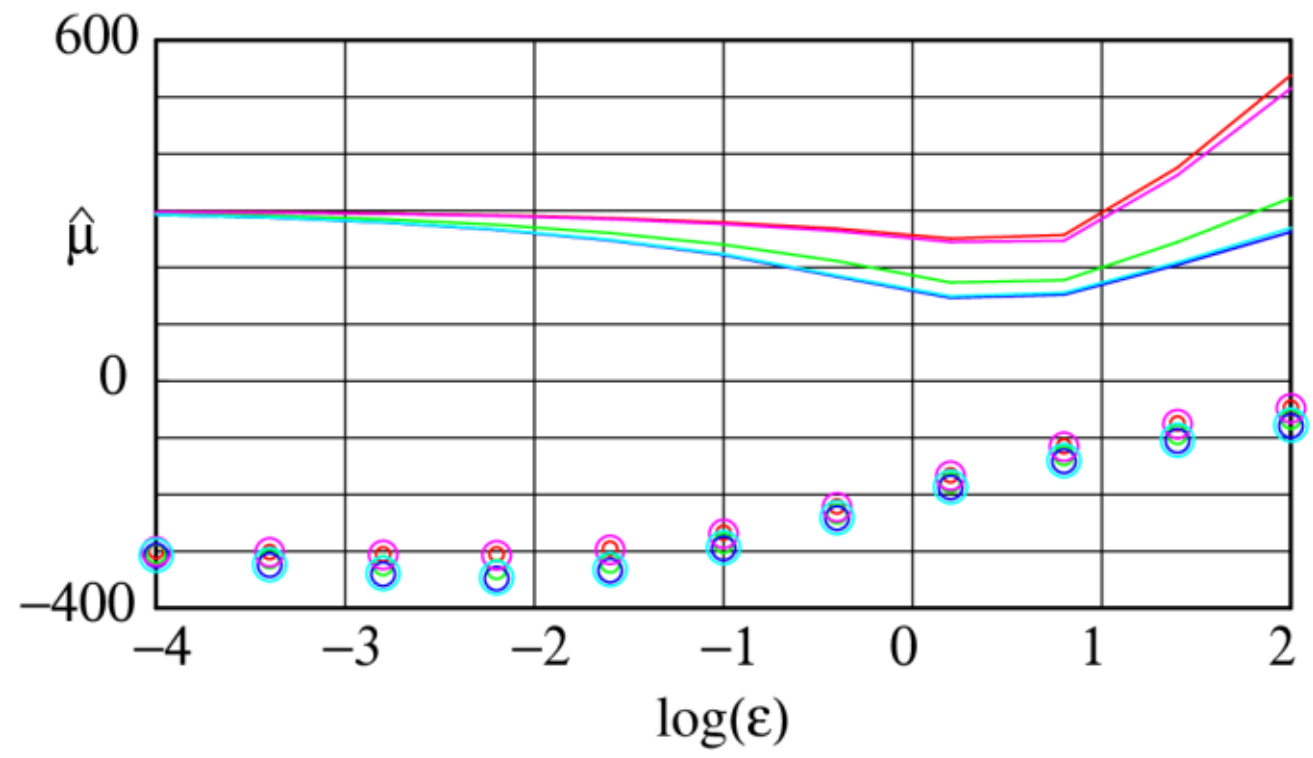

**Figure 7.** Computed $\mu = \hat{\mu}$ ($m^3/s$) values of the distorted supOU processes for different values of $\varepsilon$ in **Upper-bound** and **Lower-bound cases**. Same color legends as in **Figure 5**.

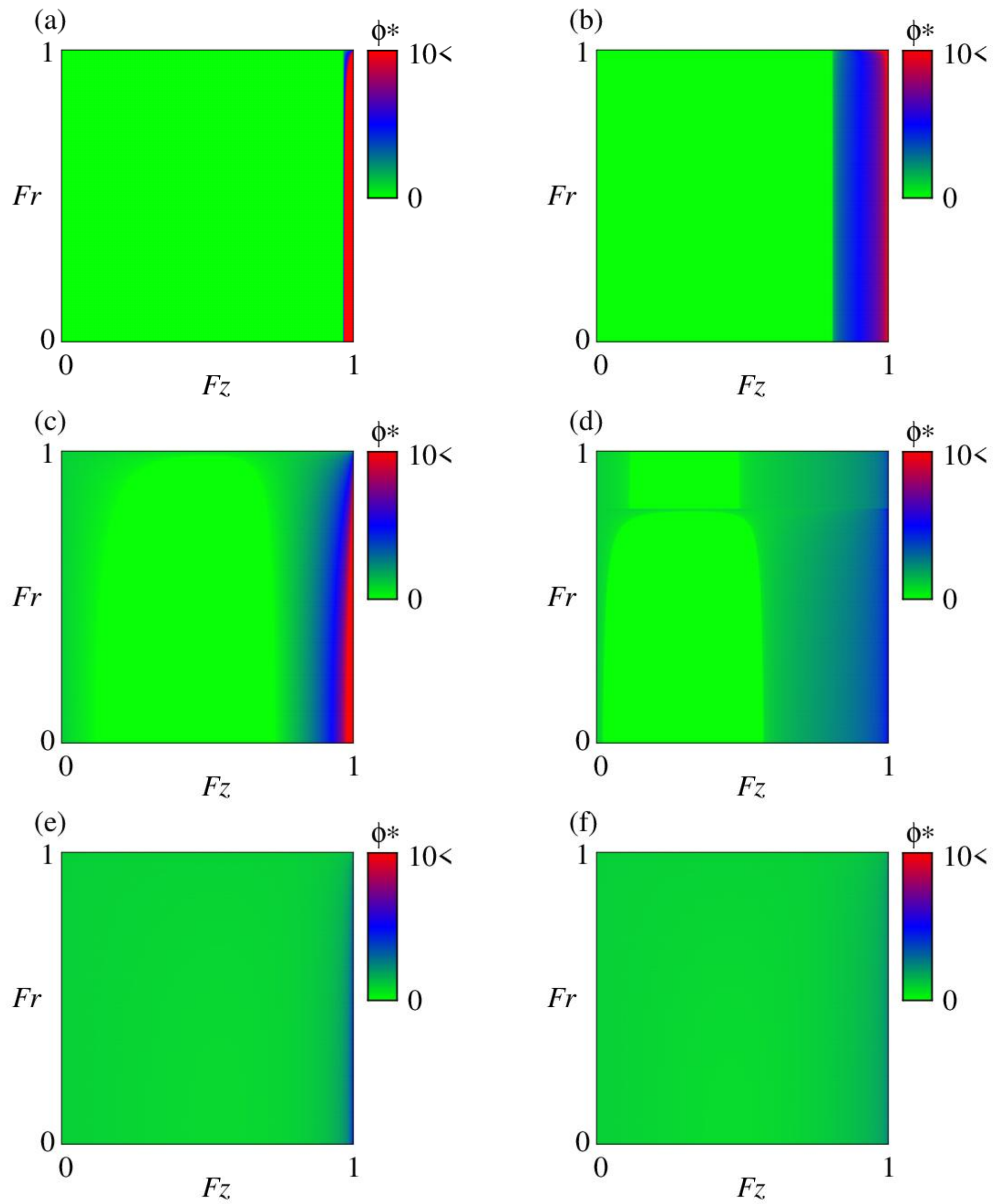

**Figure 8.** Computed $\phi^{*}$ in **Upper-bound case** with increasing and decreasing $\alpha$: (a) increasing with $\varepsilon = 100$, (b) decreasing with $\varepsilon = 100$, (c) increasing with $\varepsilon = 6.31$, (d) decreasing with $\varepsilon = 6.31$, (e) increasing with $\varepsilon = 0.1$, (f) decreasing with $\varepsilon = 0.1$.

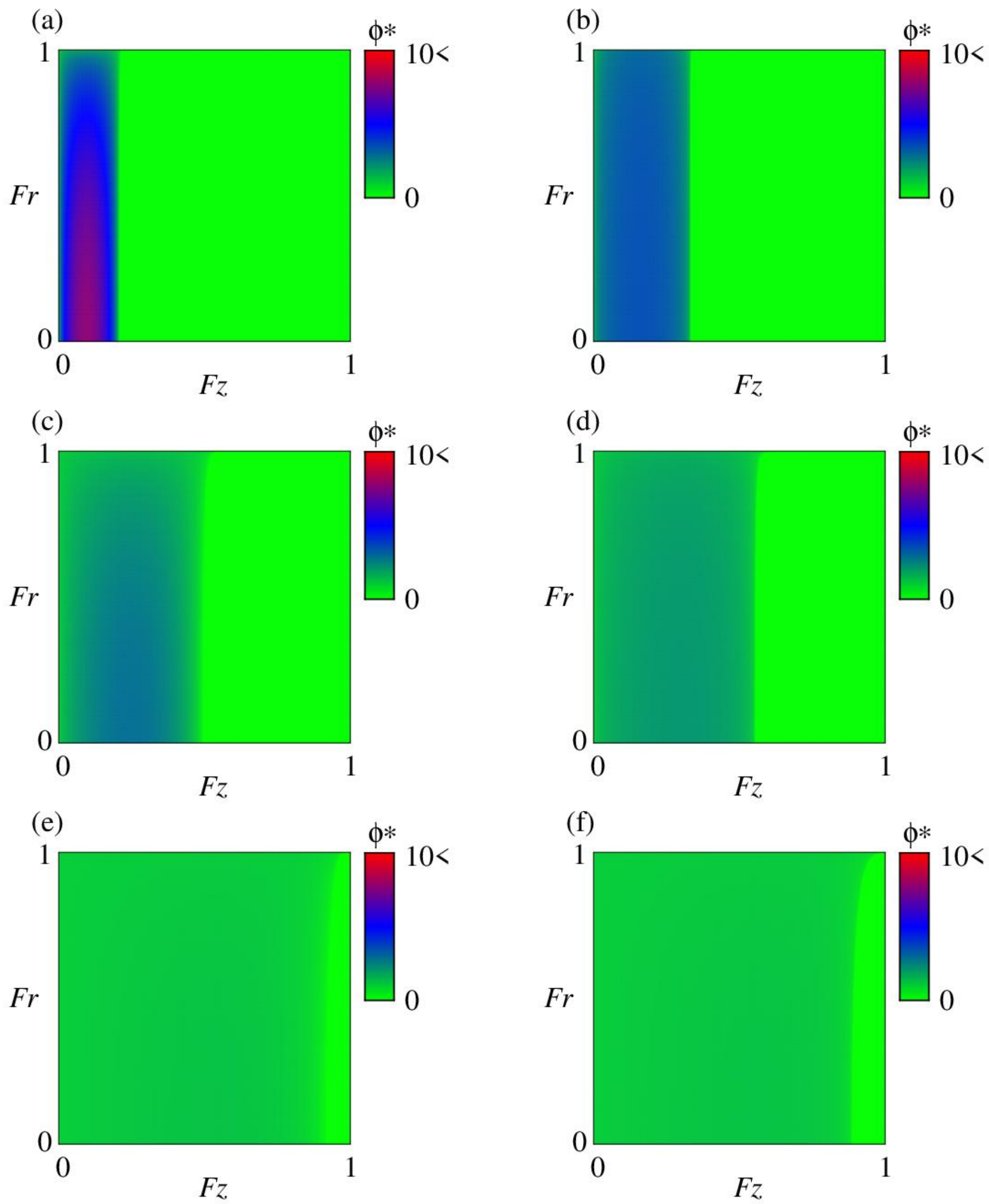

**Figure 9.** Computed $\phi^*$ in **Lower bound case** with increasing and decreasing $\alpha$ : (a) increasing with $\varepsilon = 100$, (b) decreasing with $\varepsilon = 100$, (c) increasing with $\varepsilon = 6.31$, (d) decreasing with $\varepsilon = 6.31$, (e) increasing with $\varepsilon = 0.1$, (f) decreasing with $\varepsilon = 0.1$.

#### 4.3.2 Distorted supOU processes

We computationally simulate the supOU processes for **Upper- and Lower-bound cases** by focusing on threshold crossing dynamics. Streamflow discharge itself as a continuous-state variable is important in hydrology, but its threshold statistics, such as those of high- and low-flow periods, also play roles in engineering applications such as plant uprooting [64,65], both of which are suggested to occur only during high-flow durations. Statistical analysis of flow durations is also important for freshwater fish species such as salmonids, some of which are found in the mountains reaches of the Tedori River [66]. Moreover, too little streamflow discharge is critical for water intake facilities such as those for hydropower generation [67], which is also the case around the Kazarashi station[3].

Here, a high-flow period (resp., low-flow period) is a period during which the streamflow discharge is greater than (resp., not greater than) a prescribed threshold value $X_{\mathrm{thr}} > 0$. Statistics of these periods can be computed by sampling from a sufficiently long sample path of streamflow discharge generated from a model and are inevitably affected by model uncertainties. We study the influences of model uncertainties on the statistics of high- and low-flow periods.

The averages of streamflow discharge for selected cases are presented in **Table 2** and are used in the sequel for sensitivity analysis. These cases correspond to $\varepsilon = O\left(10^{-1}\right)$, where the average values of the distorted models are approximately ±50-60% and ±20-30% of the benchmark values, which are considered within a reasonable range considering the variations in annual average values at the study site. **Figure 10** shows computed one-year sample paths for the benchmark case and **Upper- and Lower-bound cases** with $\varepsilon = 0.631$, suggesting that the overestimation of average discharge increases both minimum and maximum discharges, and vice versa. **Figure 11** shows the computed averages of the distorted supOU processes for different values of $\varepsilon$ in **Upper- and Lower-bound cases**. We focus later on the two cases $X_{\mathrm{thr}} = 20$ (m$^3$/s) (quantile level 77.9% with the benchmark model) and $X_{\mathrm{thr}} = 160$ (m$^3$/s) (quantile level 98.4% with the benchmark model); the former corresponds to the value closely linked to the threshold for hydropower intake where the water cannot be intake during low-flow periods, whereas the latter corresponds to flood events where there may be flood disturbances to the aquatic environment and ecosystems of the river.

**Figure 12** compares the computed PDFs of the streamflow discharge among the empirical and distorted supOU processes. **Figures 13-14** show the PDFs of the durations of the high- and low-flow periods, respectively. **Table 3** summarizes the average and variance for high-flow periods with $X_{\mathrm{thr}} = 20$ (m$^3$/s) computed by the benchmark and distorted model as a representative case. Similarly, **Table 4** summarizes the results for low-flow periods. **Tables 5-6** show the results with $X_{\mathrm{thr}} = 160$ (m$^3$/s) as a representative case of a larger threshold value. The tails of the PDFs of streamflow discharge become heavier as the average discharge is more strongly overestimated, as shown in **Figure 12(b)**. In contrast, as shown in

---

[3]Tedori River Project Overview https://www.hrr.mlit.go.jp/kanazawa/chisui/doc/tedori.pdf (In Japanese. Last accessed on March 4, 2026)

**Figure 12(a)**, the maximum of the PDF near the origin becomes larger as the average discharge is more strongly underestimated. In these views, the computational results demonstrate how the discharge PDF is distorted under the worst-case uncertainties. **Figures 13-14** suggest that the statistical properties of the durations of high- and low-flow periods do not change qualitatively when model uncertainties are assumed in our framework but differ quantitatively.

Because a supOU process is a non-Markovian model that is not necessarily computationally efficient because OU processes are superposed, the computational cost at least linearly scales with the degree of freedom in the $r$ direction. Model reduction focusing on threshold discharge potentially allows the modeling of high- and low-flow periods in a simpler way, which can be approximated as Markovian if they are considered to be generated by exponential distributions. Recall that the exponential distribution has a PDF $\lambda e^{-\lambda w}$ ( $w > 0$ ), where both the mean and standard deviation are $\lambda^{-1}$, and hence, the coefficient of variation is 1. In this view, **Tables 3 and 6** suggest that durations of low-flow periods can be reasonably fitted with exponential distributions, with $\lambda^{-1}$ (h) being identified as the average duration for all the computational cases, serving as a reasonable model reduction for capturing the empirical tails (**Figure 14**); however, the same seems not to be true for high-flow periods where the coefficient of variation is approximately 3 to 6, which are significantly greater than 1. This implies that the non-Markovian nature of both benchmark and distorted supOU processes appear in modeling the high-flow periods, particularly the fat tails of the PDFs of their durations, which would not be captured by exponential distributions.

A minimum model for durations of the high-flow periods with fat tails would be the following inverse-gamma distribution:

$$c'(Y)^{-\varsigma-1}\exp\left(-\frac{1}{\vartheta\varsigma}\right),\quad Y>0 \tag{66}$$

with the normalization constant $c' > 0$, shape parameter $\varsigma$ (-) and scaling parameter $\vartheta$ (day). This PDF has a tail that behaves asymptotically as $Y^{-\varsigma-1}$ for large $Y > 0$. Theoretical and empirical averages and variances are equalized to estimate the parameters $\varsigma(>2)$ and $\vartheta$ (e.g., Table 2 in Yoshioka et al. [68]).

Applying moment fitting by using the computed average and variance of high-flow periods uniquely gives the estimated values of the parameters $\varsigma$ and $\vartheta$, with the former being close to 2 (**Tables 7-8**). **Table 9** shows the percentages of the occurrence of high-flow periods for different values of $X_{\mathrm{thr}}$, which suggests how the durations of high-flow periods (and hence low-flow periods) are affected by model uncertainties. The fitted PDFs of the durations of each period are shown in **Figure 13**, suggesting a reasonable agreement between the empirical and fitted results with a slight underestimation of the empirical tails. These PDFs, which have fat tails, correspond to the long-memory nature of the underlying (distorted) supOU processes. This implies that the two-state reduced modeling of the streamflow discharge of a supOU process with long memory needs a couple of exponential and nonexponential distributions. Indeed, if a random variable $Y$ follows Eq.(66), then its inverse $Y' = Y^{-1}$, which in our context represents the reversion speed of high-flow periods, follows the gamma distribution of the shape parameter $\varsigma$ and scaling parameter $\vartheta$. This PDF has a singularity at the origin of the form $(Y')^{\varsigma-1}$. **Tables 7-8** show that the fitted values of $\varsigma$ are at most $2+O(10^{-2})$, demonstrating that durations of high-flow periods involve a wide range of timescales corresponding to the long memory nature of the supOU processes. The parameter $\varsigma$ therefore represents the degree of non-Markovian nature, and its values are almost the same across the different computational cases with different $\varepsilon$ values. This is consistent with our theoretical finding that the weight $w$ is regularized such that the memories of the distorted supOU processes are qualitatively the same. As demonstrated above, the methodology proposed in this paper can be used to analyze various aspects of the supOU processes.

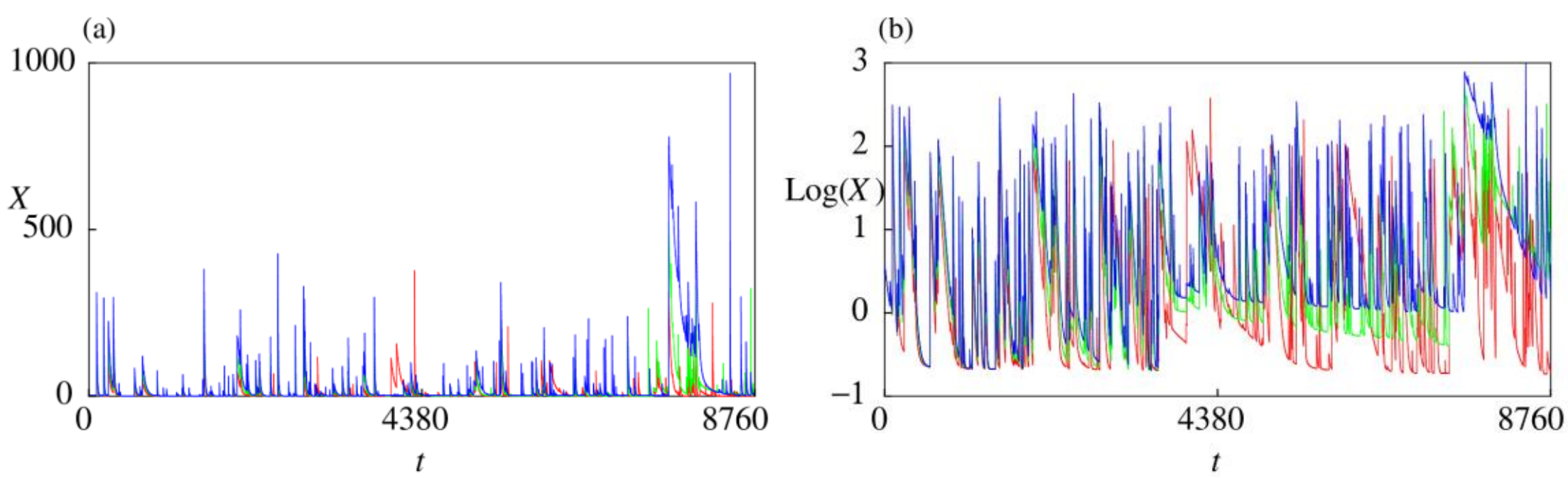

**Figure 10.** One-year ($0 \leq t \leq 8760$ (h)) sample paths of the discharge $X$ ($m^3/s$) on the (a) ordinary and (b) common logarithmic scales: the benchmark case (green), **Upper-bound case** with $\varepsilon = 0.631$ (blue), and **Lower-bound case** with $\varepsilon = 0.631$ (red).

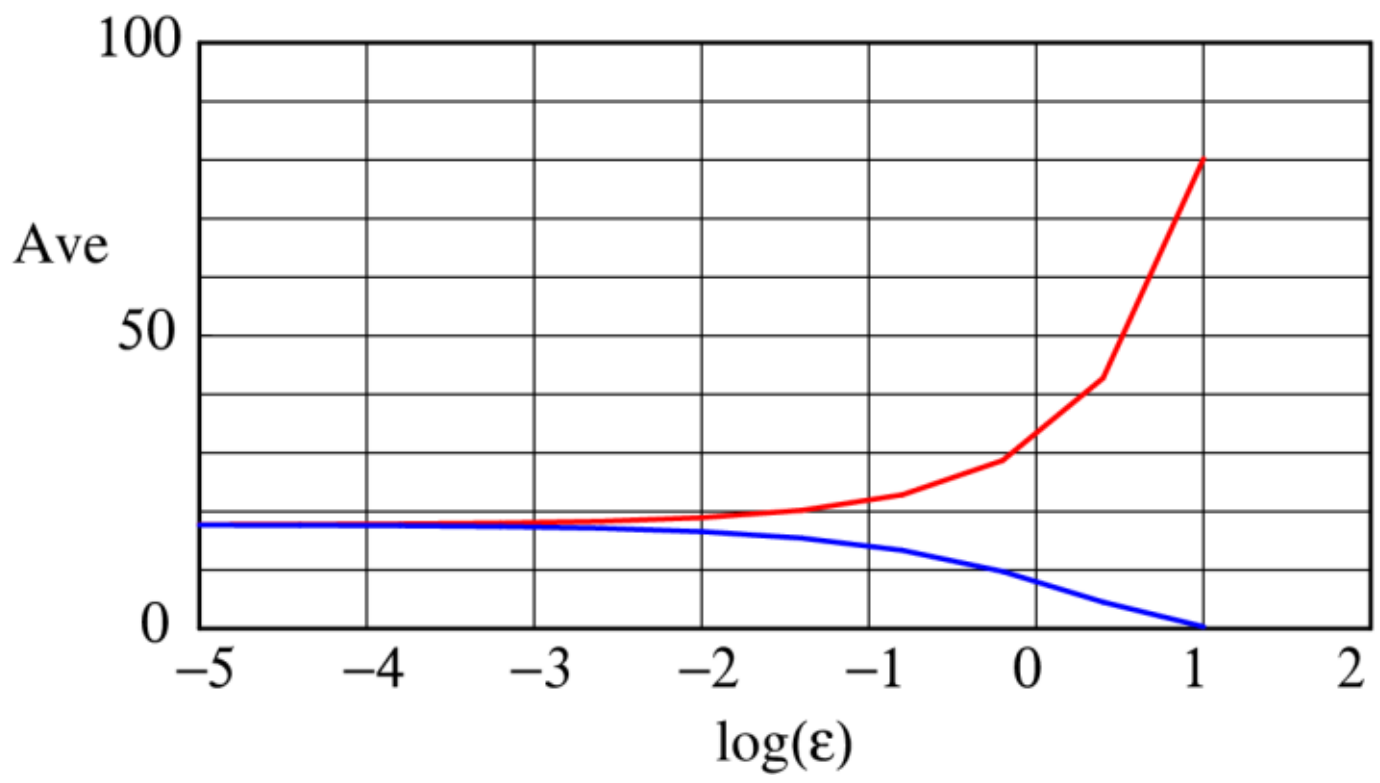

**Figure 11**. Computed average (Ave) ($m^3/s$) values of the distorted supOU processes in **Upper-bound case** (red) and **Lower-bound case** (blue).

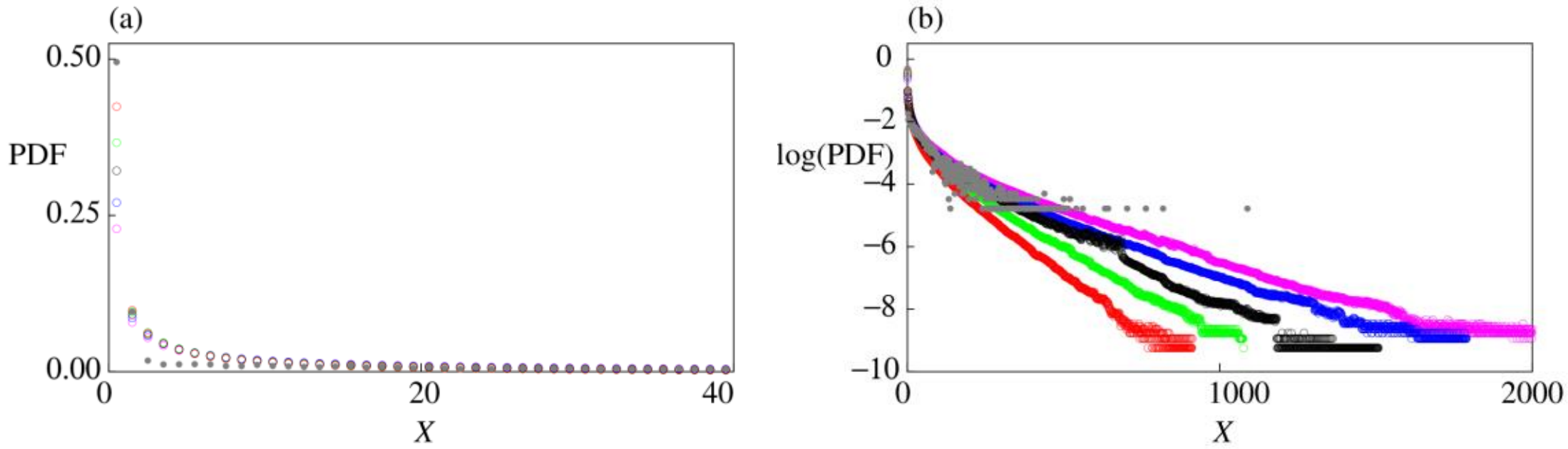

**Figure 12.** Comparison of the PDFs of streamflow discharge $X$ (m/s$^3$) among the empirical data and computational results on (a) ordinary scale and (b) common logarithmic scale: empirical data (gray), benchmark case (black), **Lower-bound case** with $\varepsilon = 0.631$ (red), **Lower-bound case** with $\varepsilon = 0.159$ (green), **Upper-bound case** with $\varepsilon = 0.159$ (blue), and **Upper-bound case** with $\varepsilon = 0.631$ (magenta).

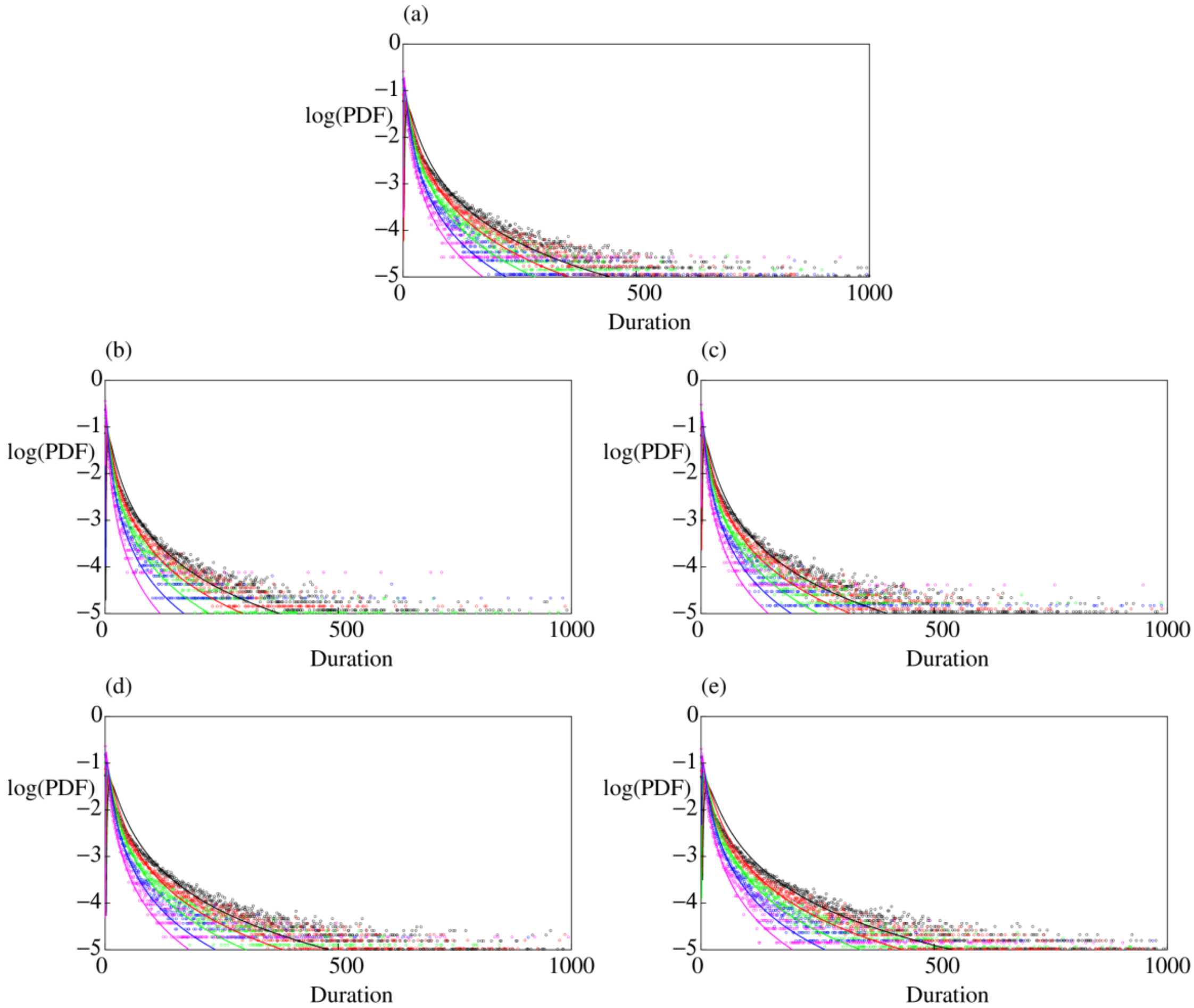

**Figure 13.** Comparison of the PDFs in the common logarithmic scale for durations (h) of the high-flow periods, where circles and curves represent empirical and fitted results, respectively: (a) Benchmark case, (b) **Lower-bound case** with $\varepsilon = 0.631$, (c) **Lower-bound case** with $\varepsilon = 0.159$, (d) **Upper-bound case** with $\varepsilon = 0.159$, (e) **Upper-bound case** with $\varepsilon = 0.631$. Color legends represent $X_{\mathrm{thr}} = 10$ (m$^3$/s) (black), $X_{\mathrm{thr}} = 20$ (m$^3$/s) (red), $X_{\mathrm{thr}} = 40$ (m$^3$/s) (green), $X_{\mathrm{thr}} = 80$ (m$^3$/s) (blue), and $X_{\mathrm{thr}} = 160$ (m$^3$/s) (magenta).

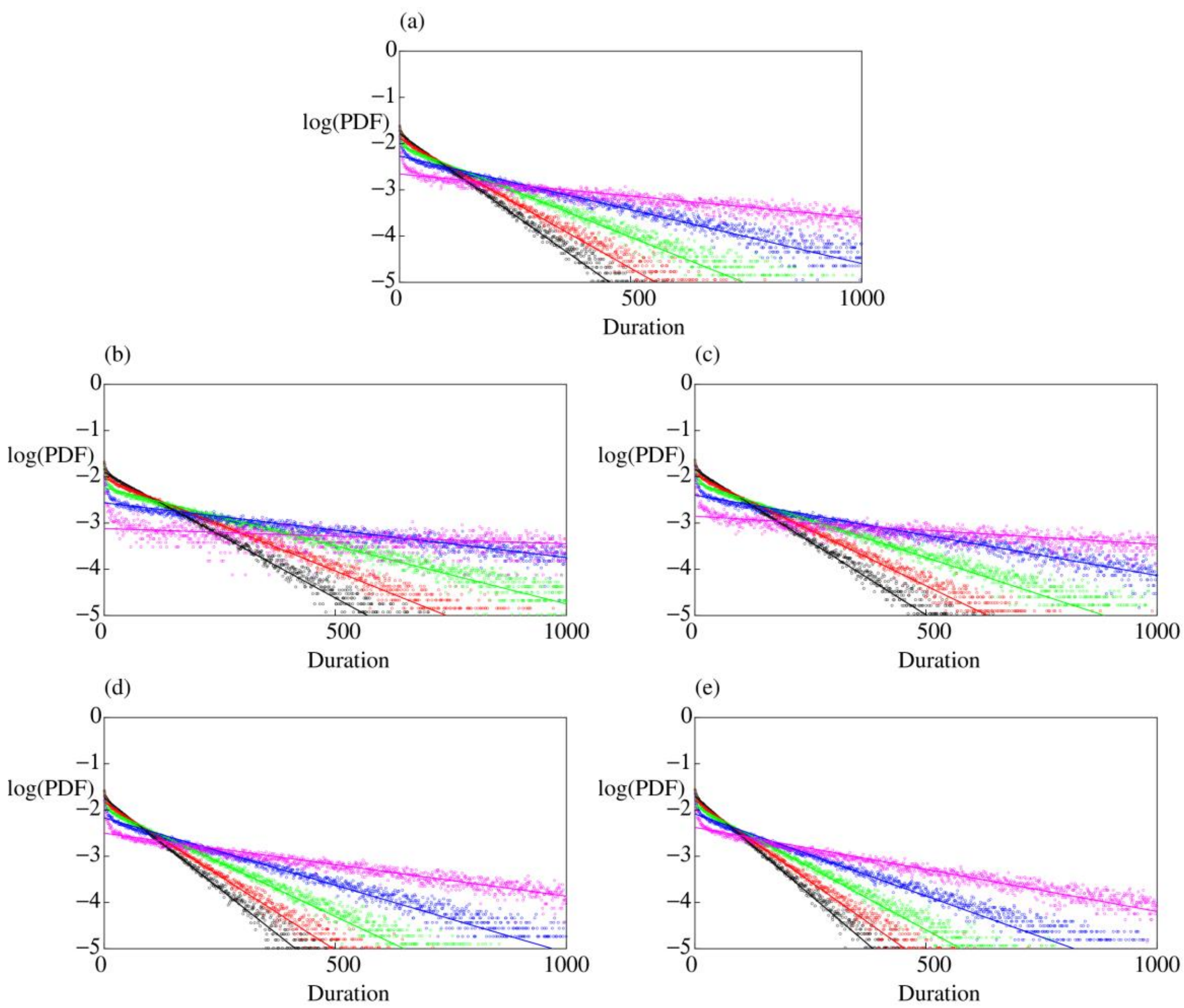

**Figure 14.** Comparison of the PDFs on the common logarithmic scale for durations of the low-flow periods, where circles and curves represent empirical and fitted results, respectively. The legends are the same as those in **Figure 13.**

**Table 2.** Average values of streamflow discharge.

| $\varepsilon$ | **Upper-bound case** | **Lower-bound case** |
|---|---|---|
| 0.159 | 22.8246 | 13.3516 |
| 0.631 | 28.6914 | 9.7366 |

**Table 3** Averages and variances of durations for high-flow periods computed by the benchmark and distorted models: $X_{\mathrm{thr}} = 20$ ($m^3/s$). "CV" indicates the coefficient of variation.

| | **Benchmark** | **Upper-bound case** | | **Lower-bound case** | |
|---|---|---|---|---|---|
| $\varepsilon$ | 0 | 0.159 | 0.631 | 0.159 | 0.631 |
| No. of samples | 176425 | 189113 | 192732 | 159940 | 140026 |
| Average (h) | 2.198.E+01 | 2.470.E+01 | 2.942.E+01 | 1.913.E+01 | 1.669.E+01 |
| Variance (h) | 1.254.E+05 | 2.778.E+04 | 7.245.E+05 | 1.028.E+04 | 1.764.E+04 |
| CV (-) | 1.611.E+01 | 6.747.E+00 | 2.894.E+01 | 5.297.E+00 | 7.959.E+00 |

**Table 4** Averages and variances of durations for low-flow periods computed by the benchmark and distorted models: $X_{\mathrm{thr}} = 20$ ($m^3/s$).

| | **Benchmark** | **Upper-bound case** | | **Lower-bound case** | |
|---|---|---|---|---|---|
| $\varepsilon$ | 0 | 0.159 | 0.631 | 0.159 | 0.631 |
| No. of samples | 176425 | 189113 | 192733 | 159941 | 140026 |
| Average (h) | 7.725.E+01 | 6.788.E+01 | 6.143.E+01 | 9.033.E+01 | 1.083.E+02 |
| Variance (h) | 6.851.E+03 | 5.230.E+03 | 4.289.E+03 | 9.302.E+03 | 1.335.E+04 |
| CV (-) | 1.071.E+00 | 1.065.E+00 | 1.066.E+00 | 1.068.E+00 | 1.066.E+00 |

**Table 5** Averages and variances of durations for high-flow periods computed by the benchmark and distorted models: $X_{\mathrm{thr}} = 160$ ($m^3/s$).

| | **Benchmark** | **Upper-bound case** | | **Lower-bound case** | |
|---|---|---|---|---|---|
| $\varepsilon$ | 0 | 0.159 | 0.631 | 0.159 | 0.631 |
| No. of samples | 37742 | 53569 | 70613 | 24301 | 13217 |
| Average (h) | 7.280.E+00 | 7.953.E+00 | 9.083.E+00 | 5.784.E+00 | 4.426.E+00 |
| Variance (h) | 1.083.E+04 | 3.054.E+03 | 2.718.E+03 | 8.676.E+02 | 2.648.E+02 |
| CV (-) | 1.430.E+01 | 6.948.E+00 | 5.739.E+00 | 5.093.E+00 | 3.676.E+00 |

**Table 6** Averages and variances of durations for low-flow periods computed by the benchmark and distorted models: $X_{\mathrm{thr}} = 160$ ($m^3/s$).

| | **Benchmark** | **Upper-bound case** | | **Lower-bound case** | |
|---|---|---|---|---|---|
| $\varepsilon$ | 0 | 0.159 | 0.631 | 0.159 | 0.631 |
| No. of samples | 37742 | 53569 | 70613 | 24301 | 13217 |
| Average (h) | 4.565.E+02 | 3.188.E+02 | 2.388.E+02 | 7.144.E+02 | 1.320.E+03 |
| Variance (h) | 2.506.E+05 | 1.219.E+05 | 7.170.E+04 | 6.016.E+05 | 1.984.E+06 |
| CV (-) | 1.097.E+00 | 1.095.E+00 | 1.121.E+00 | 1.086.E+00 | 1.067.E+00 |

**Table 7.** Fitted values of the shape parameter $\varsigma$ (-) and scaling parameter $\vartheta$ (h) for high-flow periods: $X_{\text{thr}} = 20$ (m³/s).

| | **Benchmark** | **Upper-bound case** | | **Lower-bound case** | |
|---|---|---|---|---|---|
| $\varepsilon$ | 0 | 0.159 | 0.631 | 0.159 | 0.631 |
| $\varsigma$ (-) | 2.004.E+00 | 2.022.E+00 | 2.001.E+00 | 2.036.E+00 | 2.016.E+00 |
| $\vartheta$ (h) | 2.207.E+01 | 2.525.E+01 | 2.945.E+01 | 1.982.E+01 | 1.695.E+01 |

**Table 8.** Fitted values of the shape parameter $\varsigma$ (-) and scaling parameter $\vartheta$ (h) for high-flow periods: $X_{\text{thr}} = 160$ (m³/s).

| | **Benchmark** | **Upper-bound case** | | **Lower-bound case** | |
|---|---|---|---|---|---|
| $\varepsilon$ | 0 | 0.159 | 0.631 | 0.159 | 0.631 |
| $\varsigma$ (-) | 2.005.E+00 | 2.021.E+00 | 2.030.E+00 | 2.039.E+00 | 2.074.E+00 |
| $\vartheta$ (h) | 7.315.E+00 | 8.117.E+00 | 9.359.E+00 | 6.007.E+00 | 4.754.E+00 |

**Table 9.** Percentages of high-flow periods for different $X_{\text{thr}}$ values.

| | **Benchmark** | **Upper-bound case** | | **Lower-bound case** | |
|---|---|---|---|---|---|
| $\varepsilon$ | 0 | 0.159 | 0.631 | 0.159 | 0.631 |
| $X_{\text{thr}} = 20$ (m³/s) | 22.14% | 26.66% | 32.36% | 17.47% | 13.34% |
| $X_{\text{thr}} = 160$ (m³/s) | 1.57% | 2.43% | 3.66% | 0.80% | 0.33% |

## 5. Conclusion

We proposed a novel mathematical framework for modeling the uncertainties of long-memory processes based on the divergences on Musielak–Orlicz spaces that is suited to analyzing jump-driven supOU processes. Proper Musielak–Orlicz spaces were constructed on the basis of reversion and Lévy measures. We resolved the drawback of the classical KL divergence, which fails to define cumulants subject to uncertainties well, by allowing for its state-dependence. We discussed the worst-case cumulants along with sufficient conditions for the solvability of the corresponding optimization problems. The form of the worst-case cumulants enabled us to compute them by a gradient descent method with quantization. The application study demonstrated how uncertainties in reversion and Lévy measures affect the behavior of supOU processes, particularly their PDFs and threshold statistics. This paper suggested that distorted supOU processes based on Musielak–Orlicz spaces lead to a rigorous and computable approach for modeling uncertainties in long-memory phenomena.

A limitation of this study is that we considered only supOU processes, while many other long-memory processes exist. Because the proposed mathematical framework for modeling uncertainties can simultaneously address distortions of multiple measures, it can be extended to many other jump-driven ambit and random fields driven by jump measures. Another limitation of this study is the range of application areas of supOU processes because we discussed only environmental applications; however, a wide variety of problems described by jump-driven long-memory processes potentially exist. In this view, a potential future research direction is to accumulate case studies of supOU and related processes in diverse research fields and integrate the uncertainty characteristics of each problem toward the development of a more versatile mathematical framework for modeling uncertainties in complex phenomena with long memory. For example, stability and robustness analysis of engineering systems, like those based on eigenvalue analysis [69], can be addressed from a different standpoint by using the proposed framework. Finally, because we could characterize uncertainties in the discharge data in view of Musielak–Orlicz spaces, an interesting future research topic would be mapping rivers and/or river systems via their associated Musielak–Orlicz spaces so that we can better understand aquatic system dynamics from a different perspective than before.

## Appendix

### A1. Proofs

***Proof of Proposition 1***

The proof for the first part, except for stationarity, is by a direct adaptation of Proposition 39 in Barndorff-Nielsen et al. [21] combined with the assumption that in the present case, there is no diffusion component, and the jumps are nonnegative. Stationarity follows from a direct calculation analogous to that of Proof of Theorem 3.1 in Barndorff-Nielsen [8]. The second part follows from the calculations analogous to Proof of Proposition 2.6 in Rajput and Rosinski [10].

□

***Proof of Proposition 2***

The proof relies on a Tauberian argument that connects the tail of an autocorrelation function and the singularity of the integrand (Proof of Proposition 6 in Fasen and Klüppelberg [51]) along with Eq.(10).

□

***Proof of Proposition 3***

First, we show the dual formulation in Eq.(32) by following the lines of Proof of Proposition 3.19 in Fröhlich and Williamson [36]. The difference is that we use Musielak–Orlicz spaces, whereas the literature uses classical Orlicz spaces.

The admissible set $\mathfrak{A}_{m,\varepsilon}$ is not a null set because the const function $\phi \equiv 1$ belongs to it, and hence, there exist Lagrangian multipliers $\left(\tau^{*},\mu^{*}\right) \in [0,+\infty) \times \mathbb{R}$ such that (Proof of Theorem 5.1 in Dommel and Pichler [52])

$$\overline{I}_{k,m} = \frac{mc_m}{k} \sup_{\phi \geq 0} \left\{ \mathbb{E}_m\left[z^{k-m}\phi\right] - \tau^{*}\left(\mathbb{D}'(\phi) - \frac{\varepsilon}{mc_m}\right) - \mu^{*}\left(\mathbb{E}_m[\phi] - 1\right) \right\}, \tag{67}$$

and the nonnegativity as well as lower-semi continuity of $\Phi$ shows:

$$\overline{I}_{k,m} = \frac{mc_m}{k} \inf_{\mu \in \mathbb{R},\ \tau > 0} \sup_{\phi \geq 0} \left\{ \mathbb{E}_m\left[z^{k-m}\phi\right] - \tau\left(\mathbb{D}'(\phi) - \frac{\varepsilon}{mc_m}\right) - \mu\left(\mathbb{E}_m[\phi] - 1\right) \right\}, \tag{68}$$

where the supremum is taken with respect to nonnegative random variables on $(0,+\infty)^2$. We proceed as

$$\begin{aligned} \overline{I}_{k,m} &= \frac{mc_m}{k} \inf_{\mu \in \mathbb{R},\ \tau > 0} \left\{ \frac{\varepsilon}{mc_m}\tau + \mu + \tau \sup_{\phi \geq 0}\left( \mathbb{E}_m\left[\frac{z^{k-m} - \mu}{\tau}\phi\right] - \mathbb{D}'(\phi) \right) \right\} \\ &= \frac{mc_m}{k} \inf_{\mu \in \mathbb{R},\ \tau > 0} \tau \left\{ \frac{\varepsilon}{mc_m} + \sigma + \mathbb{E}_m\left[ \Psi\left( w'(r,z), \alpha(r,z), \frac{z^{k-m}}{\tau} - \sigma \right) \right] \right\}, \end{aligned} \tag{69}$$

where we applied the replacement $\mu = \tau\sigma \in \mathbb{R}$ to obtain the last line; the second line proves Eq.(32).

Second, we show that the infimum is attained at an interior point. This does not occur if and only if $\tau$ inside the third line in Eq.(69) tends toward 0, where we use the functional similarity of the optimization problems between ours and that in the latter half of Proof of Proposition 5.5 in Dommel and Pichler [52]. If $\tau$ inside the third line in Eq.(69) tends toward 0, then $\hat{\tau}=0$ and hence

$$\begin{aligned}\bar{I}_{k,m} &= \sup_{\phi\ge 0}\left\{\mathbb{E}_m\left[z^{k-m}\phi\right]-\hat{\tau}\left(\mathbb{D}'(\phi)-\frac{\varepsilon}{mc_m}\right)-\hat{\mu}\left(\mathbb{E}_m[\phi]-1\right)\right\}\\ &= \sup_{\phi\ge 0}\left\{\mathbb{E}_m\left[z^{k-m}\phi\right]-\hat{\mu}\left(\mathbb{E}_m[\phi]-1\right)\right\}\\ &= +\infty,\end{aligned}\tag{70}$$

where the last line is due to **Assumption 2**, and we can realize it by considering $\phi$, which is concentrated at an arbitrarily large $z>0$ because $k>m$. However, we have the following upper bound of $\bar{I}_{k,m}$, leading to a contradiction:

$$\begin{aligned}\bar{I}_{k,m} &\le \frac{mc_m}{k}\tau\left\{\frac{\varepsilon}{mc_m}+\mu+\mathbb{E}_m\left[\Psi\left(w'(r,z),\alpha(r,z),\frac{z^{k-m}-\mu}{\tau}\right)\right]\right\}\Bigg|_{\tau=c>>1,\mu=0}\\ &\le \frac{mc_m}{k}c\left\{\frac{\varepsilon}{mc_m}+\mathbb{E}_m\left[\bar{\Psi}\left(w'(r,z),\alpha(r,z),z^{k-m}/c\right)\right]\right\}\\ &< +\infty,\end{aligned}\tag{71}$$

where the last line is due to $z^{k-m}\in L_{\bar{\Psi}}$ and **Assumption 3**. Therefore, the infimum is attained at an interior point $(\hat{\tau},\hat{\mu})\in(0,+\infty)\times\mathbb{R}$, which is a minimizer. For the nonnegativity of $\hat{\mu}$, because the function $\bar{\Psi}(w_0,\alpha_0,y)-y$ ($y\in\mathbb{R}$) is nonnegative, nondecreasing, and convex for any $w_0>0$ and $\alpha_0>1$, the infimum for $\mu$ in the second line of Eq.(69) is not attained for $\mu<0$ (see Proof of Theorem 4.5 in Dommel and Pichler [52]).

The maximizer $\phi^*$ of the optimization problem is obtained as a byproduct of Eq.(69). Indeed, it satisfies

$$\frac{z^{k-m}-\hat{\mu}}{\hat{\tau}}=\frac{\partial\Phi}{\partial\phi}\left(w',\alpha,\phi^*\right)\ \text{ for } (r,z) \text{ such that } \phi^*>0,\tag{72}$$

and hence $\phi^*(r,z)$ is given by Eq.(33) because $\frac{\partial\Phi}{\partial\phi}(\cdot,\cdot,\phi)$ is the inverse of $\frac{\partial\Psi}{\partial\psi}(\cdot,\cdot,\psi)$ (with respect to the third argument), and hence

$$\frac{\partial\Psi}{\partial\psi}\left(w',\alpha,\frac{z^{k-m}-\hat{\mu}}{\hat{\tau}}\right)=\frac{\partial\Psi}{\partial\psi}\left(w',\alpha,\frac{\partial\Phi}{\partial\phi}\left(w',\alpha,\phi^*\right)\right)=\phi^*.\tag{73}$$

The fact that the couple $(\hat{\tau},\hat{\mu})$ solves Eq.(34) is due to the continuous differentiability of $F$.

Finally, we show that the condition **($\mathbf{I}_k$)** with $\phi^*$ holds true. The assumption $z^{k-m}\in L_{\bar{\Psi}}$ leads to the fact that the optimizer $\phi^*$ given in Eq.(33) satisfies the condition $\mathbb{D}'(\phi^*)<+\infty$, and hence, $\phi^* \in L_{\bar{\Phi}}$. This implies that the condition **(I*k*)** with $\phi^*$ holds true because, by the Fenchel–Yoshioka inequality between $\bar{\Phi}$ and $\bar{\Psi}$ (Lemma 2.1.32 in Chlebicka et al. [33]), we obtain

$$\begin{aligned}
&\int_{r=0}^{r=+\infty}\int_{z=0}^{z=+\infty}\phi^*(r,z)z^{k-m}p_m(\mathrm{d}r\mathrm{d}z) \\
&\leq \int_{r=0}^{r=+\infty}\int_{z=0}^{z=+\infty}\bar{\Phi}\left(w',\alpha,\hat{\tau}\phi^*(r,z)\right)p_m(\mathrm{d}r\mathrm{d}z)+\int_{r=0}^{r=+\infty}\int_{z=0}^{z=+\infty}\bar{\Psi}\left(w',\alpha,z^{k-m}/\hat{\tau}\right)p_m(\mathrm{d}r\mathrm{d}z) \\
&< +\infty,
\end{aligned} \tag{74}$$

where the last line comes from $\phi^* \in L_{\bar{\Phi}}$ and $z^{k-m} \in L_{\bar{\Psi}}$ along with **Assumption 3**.

□

***Proof of Proposition 4***

Most parts of the proof here are the same as those for **Proposition 3** except for the point to prove $\hat{\mu} \in (-\infty, 0]$. This is due to $0 \leq \dfrac{\partial \Psi}{\partial \psi}(w_0,\alpha_0,\psi) \leq 1$ for any $w_0 > 0$, $\alpha_0 > 1$, and $\psi \leq 0$. Indeed, because we need $\mathbb{E}_m\left[\phi^*\right]=1$, $\hat{\mu}$ must be positive considering Eq.(37) and $-\dfrac{z^{k-m}}{\hat{\tau}} \leq 0$; otherwise, we face $\phi^* < 1$ and the condition $\mathbb{E}_m\left[\phi^*\right]=1$ is never satisfied.

□

***Proof of Proposition 5***

By **Definition 1**, if $\phi \in L_{\bar{\Phi}(w_2',\alpha_2,\cdot)}$, then $\phi \in L_{\bar{\Phi}(w_1',\alpha_1,\cdot)}$; therefore, the first part of the proof holds true. The second part is due to the following inequality: for $\phi \in L_{\bar{\Phi}(w_2',\alpha_2,\cdot)}$,

$$\|\phi\|_{(w_2',\alpha_2,\cdot)} = \inf\left\{c>0 \middle| \mathbb{E}_m\left[\bar{\Phi}\left(w_2',\alpha_2,\frac{|\phi|}{c}\right)\right]\leq 1\right\} \geq \inf\left\{c\mathbb{E}_m\left[\bar{\Phi}\left(w_1',\alpha_1,\frac{|\phi|}{c}\right)\right]\leq 1\right\} = \|\phi\|_{(w_1',\alpha_1,\cdot)}. \tag{75}$$

□

**A2. Numerical method**

We apply the following gradient descent with momentum to solve the system (34), which is a version of the Nesterov method (e.g., Walkington [70]):

$$\begin{aligned}
\lambda_{n+1} &= \omega\lambda_n - \eta_n\frac{\partial F\left(\tau_n+\omega\lambda_n,\mu_n+\omega\rho_n\right)}{\partial\tau}, \quad n=0,1,2,\ldots \\
\tau_{n+1} &= \tau_n + \lambda_n
\end{aligned} \tag{76}$$

and

$$\begin{aligned}
\rho_{n+1} &= \omega\rho_n - \kappa_n\frac{\partial F\left(\tau_n+\omega\lambda_n,\mu_n+\omega\rho_n\right)}{\partial\mu}, \quad n=0,1,2,\ldots \\
\mu_{n+1} &= \mu_n + \rho_n
\end{aligned} \tag{77}$$

with an initial guess $(\tau_0,\mu_0)\in(0,+\infty)\times\mathbb{R}$ and $(\lambda_0,\rho_0)\in\mathbb{R}^2$. Here, the subscript $n$ represents the time step, $\omega(=0.95)$ is a relaxation factor, $\lambda$ is momentum for numerical $\tau$, $\rho$ is momentum for numerical $\mu$, and $\eta$ and $\kappa$ are learning rates for $\tau$ and $\mu$, respectively. The partial derivatives are discretized by using Eq.(58). The recursions (76) and (77) are terminated if

$$\left|\frac{\partial F(\tau_n+\omega\lambda_n,\mu_n+\omega\rho_n)}{\partial\mu}\right|,\left|\frac{\partial F(\tau_n+\omega\lambda_n,\mu_n+\omega\rho_n)}{\partial\tau}\right|\le\delta, \tag{78}$$

where $\delta>0$ is a prescribed error threshold and $(\tau_{n+1},\mu_{n+1})$ is considered an approximation of $(\hat{\tau},\hat{\mu})$. We empirically found that the following learning rates work well for our computational case:

$$\eta_n=2.5\varepsilon^{-1} \text{ and } \kappa_n=0.25\tau_n^{-1} \tag{79}$$

based on the preliminary finding that $\tau_n$ is large, such as $O(10^2)$ or larger, and a small $\eta_n$ is needed for the stability of the algorithm when $\varepsilon$ is large.

Here, we report the computational costs of the proposed numerical algorithm, where we set $M=512$. **Tables A1-A2** show the iteration count and computed $(\hat{\tau},\hat{\mu})$ for different values of $\varepsilon$ for $\alpha\equiv 2.5,\ 4.0,\ 5.0$. Software (C++ codes) to obtain these results have been uploaded at the following repository: https://github.com/HidekazuYoshioka/CSF_Musielak.git. Here, for each $\alpha$, we start the computation from $\varepsilon=10^{2-0.6\times0}$ and the initial guess at $\varepsilon=10^{2-0.6\times1}$ is that of the numerical solution at $\varepsilon=10^{2-0.6\times0}$. This procedure is iterated for subsequent computations for $\varepsilon=10^{2-0.6i}$ ($i=2,3,\ldots$). The initial guess at $\varepsilon=10^{2-0.6\times0}$ is $(\tau,\mu)=(500,100)$ for both **Upper- and Lower-bound cases**. The iteration count needed to obtain a numerical solution increases as $\varepsilon$ decreases under the present computational setting, where we should use a large learning parameter due to a small $\hat{\tau}$, as indicated in **Section 3.3**. This implies the need for a larger number of iterations with (sufficiently small) constant $\eta$ and $\kappa$ values. The iteration counts among different $\alpha$ values are not critically different, and the computational costs for **Upper- and Lower-bound cases** are comparable. Complexity of the numerical method is approximately given by $\log_{10}(\text{Iteration})=c_1-c_2\log_{10}\varepsilon$ with some $c_1,c_2>0$ for each $\alpha$ value, as reported in **Tables A1-A2**. **Table A3** compares the exact and computed variances for different values of $M$. The exact variance is directly obtained from Eq.(7). The relative error at $M=512$ used in the main text is 0.056%, and the convergence order of the discretization is approximately 1; the relative error is $O(M^{-1})$ because the error is almost halved by doubling $M$.

**Table A1.** Iteration counts for convergence along with the fitting results: **Upper-bound case**.

| $i$ | $\varepsilon$ | $\alpha \equiv 2.5$ | $\alpha \equiv 4.0$ | $\alpha \equiv 5.0$ |
|---|---|---|---|---|
| 0 | 1.000.E+02 | 727 | 611 | 621 |
| 1 | 2.512.E+01 | 632 | 548 | 517 |
| 2 | 6.310.E+00 | 1092 | 569 | 573 |
| 3 | 1.585.E+00 | 1738 | 1467 | 1263 |
| 4 | 3.981.E-01 | 3918 | 3378 | 3154 |
| 5 | 1.000.E-01 | 8226 | 7645 | 7242 |
| 6 | 2.512.E-02 | 16660 | 16435 | 15917 |
| 7 | 6.310.E-03 | 32912 | 34033 | 33628 |
| 8 | 1.585.E-03 | 63803 | 68329 | 68640 |
| 9 | 3.981.E-04 | 121944 | 133869 | 136151 |
| 10 | 1.000.E-04 | 230607 | 257622 | 264250 |
| | $c_1$ | 3.49.E+00 | 3.42.E+00 | 3.40.E+00 |
| | $c_2$ | 4.57.E-01 | 4.90.E-01 | 4.94.E-01 |
| | $R^2$ | 9.85.E-01 | 9.76.E-01 | 9.73.E-01 |

**Table A2.** Iteration counts for convergence along with the fitting results: **Lower-bound case**.

| $i$ | $\varepsilon$ | $\alpha \equiv 2.5$ | $\alpha \equiv 4.0$ | $\alpha \equiv 5.0$ |
|---|---|---|---|---|
| 0 | 1.000.E+02 | 539 | 522 | 519 |
| 1 | 2.512.E+01 | 501 | 501 | 501 |
| 2 | 6.310.E+00 | 521 | 501 | 501 |
| 3 | 1.585.E+00 | 1004 | 950 | 944 |
| 4 | 3.981.E-01 | 3570 | 3427 | 3301 |
| 5 | 1.000.E-01 | 9445 | 10003 | 9832 |
| 6 | 2.512.E-02 | 20910 | 22434 | 25415 |
| 7 | 6.310.E-03 | 38327 | 29104 | 13366 |
| 8 | 1.585.E-03 | 68744 | 82046 | 97235 |
| 9 | 3.981.E-04 | 126670 | 106563 | 181052 |
| 10 | 1.000.E-04 | 235436 | 267078 | 273483 |
| | $c_1$ | 3.40.E+00 | 3.39.E+00 | 3.38.E+00 |
| | $c_2$ | 5.02.E-01 | 5.04.E-01 | 5.14.E-01 |
| | $R^2$ | 9.70.E-01 | 9.66.E-01 | 9.50.E-01 |

**Table A3.** Maximum errors between exact and computed variances. The exact variance is 1472.48 ($m^6/s^2$).

| Resolution $M$ | Computed variance ($m^6/s^2$) | Relative error (%) |
|---|---|---|
| 32 | 1459.32 | 0.894 |
| 64 | 1465.89 | 0.448 |
| 128 | 1469.19 | 0.223 |
| 256 | 1470.83 | 0.112 |
| 512 | 1471.66 | 0.056 |

### A3. Influences of regularization

**Figure A1** shows the computed variance values of the distorted supOU processes for different values of $\varepsilon$ in **Upper- and Lower-bound cases** with different regularization parameter $W$ values, demonstrating that the regularization does not quantitatively affect the variance values under uncertainties. We empirically found that specifying a weaker regularization (i.e., larger $W$ value) results in slower convergence or nonconvergence of the computational method, which is due to the coexistence of terms ranging across more than several tens of orders in Eq.(58). Therefore, the regularization avoids this computational failure. The development of a stable computational methodology to address integrands without regularizations is currently ongoing.

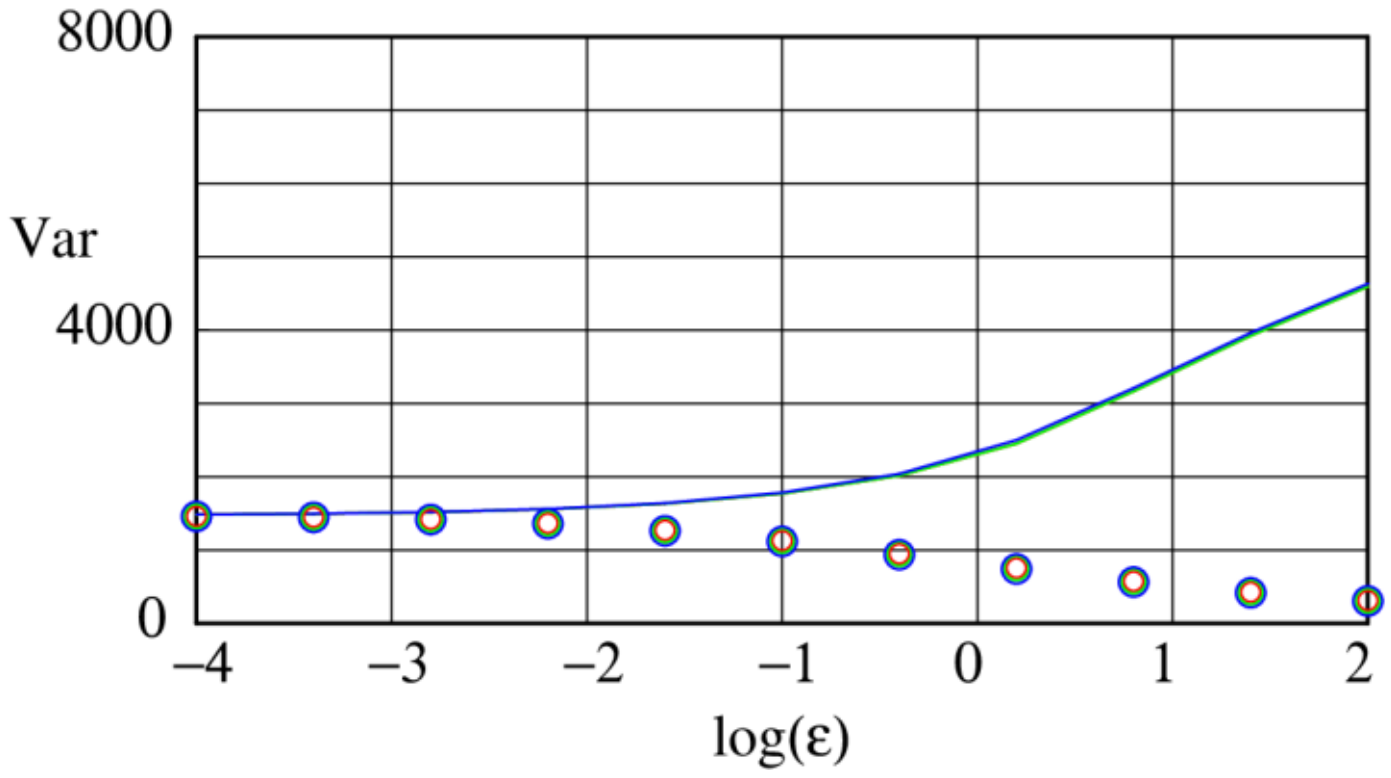

**Figure A1.** Computed variance (Var) ($m^6/s^2$) values of the distorted supOU processes in **Upper-bound case** (curves) and **Lower-bound case** (circles) with different regularization parameter $W$ (h) values: $W = 0.1$ (red), $W = 1$ (green), and $W = 10$ (blue).